\documentclass[11pt]{article}
\usepackage{amssymb,amsmath,latexsym}

\usepackage{epsfig}

\newtheorem{thm}{Theorem}[section]
\newtheorem{defn}[thm]{Definition}
\newtheorem{corollary}[thm]{Corollary}

\newtheorem{lemma}[thm]{Lemma}

\newtheorem{remark}[thm]{Remark}
\newtheorem{example}[thm]{Example}
\newtheorem{assumption}[thm]{Assumption}

\newcommand{\pf}{\noindent{\bf Proof.} }
\def\qed{{\hfill $\Box$ \bigskip}}

\newcommand\cbrk{\text{$]$\kern-.15em$]$}}
\newcommand\opar{\text{\,\raise.2ex\hbox{${\scriptstyle
|}$}\kern-.34em$($}}
\newcommand\cpar{\text{$)$\kern-.34em\raise.2ex\hbox{${\scriptstyle |}$}}\,}

\def\wh{\widehat}
\def\wt{\widetilde}
\def\<{\langle}
\def\>{\rangle}
\def\eps{\varepsilon}
\def\E{{\mathbb E}}

\newcommand\bL{\mathbb{L}}
\newcommand\bR{\mathbb{R}}
\newcommand\bH{\mathbb{H}}

\newcommand\bD{\mathbb{D}}
\newcommand\bS{\mathbb{S}}

\newcommand\cA{\mathcal{A}}
\newcommand\cB{\mathcal{B}}

\newcommand\cF{\mathcal{F}}
\newcommand\cH{\mathcal{H}}

\newcommand\cP{\mathcal{P}}
\newcommand\cR{\mathcal{R}}

\newcommand\cZ{\mathcal{Z}}
\def\loc{{\rm loc}}

\def\wh{\widehat}
\def\wt{\widetilde}

\newcommand{\mysection}[1]{\section{#1}
\setcounter{equation}{0}}

\begin{document}

\title{\bf  An
$L_p$-theory of stochastic parabolic equations   with the random fractional Laplacian
driven by L\'evy processes}

\date{}

\author{
Kyeong-Hun Kim\footnote{Department of
Mathematics, Korea University, 1 Anam-dong, Sungbuk-gu, Seoul 136-701, Republic of Korea. \,\,E-mail:  \texttt{kyeonghun@korea.ac.kr}. The research of this
author  was supported by Basic Science Research Program through the National Research Foundation
of Korea(NRF) funded by the Ministry of Education, Science and Technology (20090087117).}
\qquad \hbox{\rm and} \qquad
Panki Kim\footnote{
Department of Mathematical Sciences and Research Institute of Mathematics,
Seoul National University,
San56-1 Shinrim-dong Kwanak-gu,
Seoul 151-747, Republic of Korea.
\,\,
E-mail: \texttt{pkim@snu.ac.kr}.
This work was supported by Basic
Science Research Program through the National Research Foundation of
Korea(NRF) grant funded by the Korea
government(MEST)(2010-0001984).      }
}

\maketitle

\begin{abstract}
In this paper we give an $L_p$-theory for stochastic parabolic equations with random fractional Laplacian operator. The driving noises are general L\'evy processes.
\end{abstract}

\vspace*{.125in}


\vspace*{.125in}


\mysection{Introduction}

Let $d, m \ge 1$ be positive integers,  $p\in [2,\infty)$ and $\alpha \in (0,2)$. As usual $\bR^{d}$ stands for the Euclidean space of points
$x=(x^{1},...,x^{d})$.
We will use $dx$  to
denote the Lebesgue measure in either $\bR^d$ or $\bR^m$, which  is clear in each context.

 In this article we are
dealing with
 $L_p$-theory
 of the stochastic partial differential equations of the type
\begin{equation}   \label{eqn 6.15.1}
du=\left(a(\omega,t)\Delta^{\alpha/2}u+f(u)\right)dt\,
+\sum_{k=1}^{\infty}g^k(u) \cdot dZ^k_t, \quad u(0)=u_0
\end{equation}
given for $\omega\in \Omega, t\geq0$  and $x\in \bR^d$.  Here  $\Omega$ is a probability space,
$\Delta^{\alpha/2}$ is the fractional Laplacian defined  in \eqref{e:fl},  $Z^k_t$ are independent
$m$-dimensional L\'evy processes,   and
 the functions $f$ and  vector-valued function $g^k=(g^{k,1}, \dots g^{k,m})$
depend on $(\omega,t,x,u)$ satisfying  certain continuity conditions.
Our result will cover the case
\begin{align*}
f(u)&=b(\omega,t,x)\Delta^{\beta_1/2}u+c^i(\omega,t,x)u_{x^i}I_{\alpha>1}+d(\omega,t,x)u+f_0,\\
 g^{k,j}(u)&=\sigma^{k,j}(\omega,t,x)\Delta^{
 \beta^j_2/2}u+\nu^{k,j}(\omega,t,x)u+g^{k,j}_0, \quad j=1, \dots m.
 \end{align*}
where $\beta_1<\alpha$ and $\beta^j_2<\alpha/2$ (see Assumptions \ref{ass nonlinear} and \ref{ass nonlinear1}).

An $L_p$-theory of \eqref{eqn 6.15.1} is introduced in \cite{CL} for the case that  $a(\omega,t)=1$ and  are only finitely many Winer processes  appear in the equation.     The approach in \cite{CL} cannot cover the case when there are infinitely many Wiener processes, and the assumptions on $g$ in \cite{CL} are stronger than
 conditions in our paper (See Remark \ref{remark kijung} below). Moreover  equations driven by jump processes are not considered in \cite{CL}.  A  H\"older space theory for more general (but non-random) integro-differential equations driven by Hilbert space-valued Wiener process is given   in \cite{M2} (also see \cite{M1} for a deterministic equation).
  Even though the main result in  \cite{M2}
 provide a nice H\"older regularity of the solution to such problem,
due to the H\"older-type function spaces defined there, assumptions on $f$ and $g$ are quite strong.
Furthermore in \cite{M2} the equations with discontinuous L\'evy processes are not considered.
We emphasize that
the approach of this paper, based on $L_p$ theory in \cite{Kr99},  is different from \cite{M2}.
 Our results include the case when  $f$ and $g$ are only distributions and the number of derivatives of $f$ and $g$ are negative and fractional.
On the other hand if $f$ and $g$ are sufficiently smooth in $x$ then  Sobolev embedding theorem combined with our $L_p$-theory  gives pointwise  H\"older continuity of the solution even when
$Z^k$ are general L\'evy processes.

$L_p$-theory for second-order stochastic parabolic equations  driven by Wiener processes was first established by Krylov \cite{Kr99}. Recently in \cite{C-K2} $L_p$ regularity theory for second-order stochastic parabolic equations  driven by L\'evy processes is discussed.

In this paper, we establish an $L_p$-theory for stochastic parabolic equations  with the random fractional Laplacian
driven by arbitrary L\'evy processes.  Our result includes the case when the equation is driven by L\'evy space-time white noise (see Theorems \ref{thm last1} and \ref{thm last2}).
Among main tools used in the article to study $L_p$-regularity theory are Burkholder-Davis-Gundy inequality and  a
parabolic version of Littlewood-Paley inequality  for the fractional Laplacian operator introduced   \cite{Ildoo}.

The organization of this article is as follows. First, in section \ref{section 2},   we prove  uniqueness and existence  results  of equation \eqref{eqn 6.15.1} driven by Wiener processes in the space $L_p(\Omega\times [0,T],H^{\gamma+\alpha/2}_p)$ (Theorem \ref{e:twcase}). Here $p\in [2,\infty)$ and $\gamma \in \bR$.  In section \ref{section 3} we extend Theorem \ref{e:twcase}
for the case when   $Z^k_t$  are L\'evy processes and $Z^k_t$ have finite $p$-th moments (see condition (\ref{main assump})).
In section \ref{section 4}, the condition (\ref{main assump}) is weakened, and the uniqueness and existence results are proved in the space $L_{p,\text{loc}}(\Omega\times [0,T],H^{\gamma+\alpha/2}_p)$. The
 condition (\ref{main assump}) can be completely dropped if  only finitely many L\'evy processes
appear in the equation.

  If we write $c=c(...)$,
this means that the constant $c$ depends only on what are in
parenthesis. The constant $c$  stands for constants whose
values are unimportant and which may change from one appearance to
another.
The dependence of the lower case constants on the  dimensions $d, m$ may not be mentioned explicitly.
We will use ``$:=$" to denote a definition, which is
read as ``is defined to be". For $a, b\in \bR$, $a\wedge b:=\min
\{a, b\}$ and $a\vee b:=\max\{a, b\}$.
Let $C^{\infty}_0(\bR^d)$ be the collection of
smooth
functions with compact supports in $\bR^d$.
Most of functions we discuss in this paper are random (depend on $\omega \in \Omega$).
For notational convenience, we suppress the dependency on $\omega$ in most of expressions

\mysection{Stochastic Parabolic equations   with the random fractional Laplacian driven by Wiener processes}
                          \label{section 2}

Let $(\Omega,\cF,P)$ be a complete probability space,
$\{\cF_{t},t\geq0\}$ be an increasing filtration of $\sigma$-fields
$\cF_{t}\subset\cF$, each of which contains all $(\cF,P)$-null sets.
We assume that on $\Omega$ we are given independent  one-dimensional
Wiener processes
 $W^{1}_{t},W^{2}_{t},...$
relative to $\{\cF_{t},t\geq0\}$. Let $\cP$ be the predictable
$\sigma$-field generated by $\{\cF_{t},t\geq0\}$.

Let  $p(t,x)$, where $t>0$, denote the
inverse Fourier
 transform of  $e^{- |\xi|^\alpha t}$ in $\bR^d$, that is,
$$
  p(t,x):= \frac{1}{(2\pi)^{d/2}}\int_{\bR^d} e^{i \xi \cdot x }e^{-
|\xi|^\alpha t} d\xi.
$$
 For a suitable function $g$ and $t>0$,  define the corresponding convolution operator
\begin{equation}
                            \label{06.09.1}
T_tg(x) := (p(t,\cdot) * g(\cdot))(x):=\int_{\bR^d} p(t,x-y)g(y)dy,
\end{equation}
and define
\begin{equation}\label{e:fl}
\partial^{\alpha}_xg(x)={\Delta}^{\frac{\alpha}{2}}g(x) = -(-{\Delta})^{\frac{\alpha}{2}}g(x):=
\cF^{-1}(-|\xi|^{\alpha} \cF(g)(\xi))(x),
\end{equation}
where
$\cF(g)(\xi)=\hat{g}(\xi):=\frac{1}{(2\pi)^{d/2}}\int_{\bR^d}e^{-i\xi\cdot
x}g(x)dx$ is the Fourier transform of $g$ in $\bR^d$.

In this section we study the nonlinear equations of the type
\begin{equation}   \label{eqn 3.17.1}
du=\left(a(\omega,t)\Delta^{\alpha/2}u+f(u)\right)\,dt
+\sum_{k=1}^{\infty}g^k(u)dW^k_t, \quad u(0)=u_0,
\end{equation}
where $a(\omega,t)\in (\delta,\delta^{-1})$ for some $\delta>0$, and $f(u)=f(\omega,t,x,u)$ and $g^k(u)=g^k(\omega,t,x,u)$ satisfy certain 
continuity conditions, which we will put below.

First  we introduce some stochastic Banach spaces.
Let $(\phi, \psi):=\int_{\bR^d} \phi(x) \psi(x)dx$ and for $p \ge 1$,
$$
L_p =L_p(\bR^d):=\{ \phi : \bR^d \to \bR,  \| \phi \|^p_p:=\int_{\bR^d} |\phi(x)|^p dx < \infty   \}.$$
 For $n=0,1,2,...$,
define
$$
H^n_p=H^n_p(\bR^d):=\left\{u: u, Du,...,D^n u\in L_p (\bR^d)
\right\}.
$$
In general, for $\gamma \in \bR$ define the space
$H^{\gamma}_p=H^{\gamma}_p(\bR^d)=(1-\Delta)^{-\gamma/2}L_p$ (called
the space of Bessel potentials or the Sobolev space with fractional
derivatives) as the set of all distributions $u$ on $\bR^d$ such
that $(1-\Delta)^{\gamma/2}u\in L_p$.
For $u\in H^\gamma_p$, we define
\begin{equation}   \label{e:2.2}
\|u\|_{H^{\gamma}_p}:=\|(1-\Delta)^{\gamma/2}u\|_p
:=\|\cF^{-1}[(1+|\xi|^2)^{\gamma/2}\cF(u)(\xi)]\|_p,
\end{equation}
where $\cF$ is the Fourier transform in $\bR^d$.
For $\ell_2$-valued $g=(g^1, g^2, \dots)$,
we define
$$\|g\|_{H^{\gamma}_p(\ell_2)}:=\|(1-\Delta)^{\gamma/2}g|_{\ell_2}\|_p
:=\|\cF^{-1}[(1+|\xi|^2)^{\gamma/2}\cF(g)(\xi)]|_{\ell_2}\|_p.
$$

Let $\overline{\cP}$ be the
completion of $\cP$ with respect to $dP\times dt$, and
  $\bH^{\gamma}_p(T):=L_p(\Omega\times [0,T],\overline{\cP},H^{\gamma}_p)$, that is, $\bH^{\gamma}_p(T)
  $ is the set of   all $\overline{\cP}$-measurable processes
$u:
\Omega \times [0,T]
\to H^{\gamma}_p$ so that
$$
\|u\|_{\bH^{\gamma}_p(T)}:= \left( \E \left[ \int^T_0
\,\|u(\omega, t)\|^p_{H^{\gamma}_p}\,dt \right] \right)^{1/p}<\infty.
$$

\begin{lemma}
                         \label{multiplier}
   For any $\beta>0$, $\eta^1_{\beta}(\xi):=\frac{(1+|\xi|^2)^{\beta/2}}{1+|\xi|^{\beta}}$, $\eta^2_{\beta}=(\eta^1_{\beta})^{-1}$,  $\eta^3:=\frac{|\xi|^{\beta}}{1+|\xi|^{\beta}}$ and $\eta^4:=\frac{|\xi|^{\beta}}{(1+|\xi|^{2})^{\beta/2}}$
   are $L^p(\bR^d)$-multipliers, that is,
   $$
   \|\cF^{-1}\left( \eta^i_{\beta}(\xi) (\cF u)(\xi)\right)\|_{L_p}\leq c(p,\beta)\|u\|_{L_p}, \quad \quad i=1,2,3,4.
   $$
   \end{lemma}
     \pf
     See Theorem 0.2.6 of \cite{Sogge} (also see the remark below the theorem).
     \qed

Lemma \ref{multiplier} easily yields the following results.
\begin{corollary}
                            \label{cor multiplier}
(i) Let $\gamma\geq 0$. There exists a constant $c=c(\gamma)>0$ so that
$$
c\|u\|_{H^{\gamma}_p}\leq (\|u\|_{L_p}+\|\partial^{\gamma/2}_xu\|_{L_p})\leq c^{-1}\|u\|_{H^{\gamma}_p}.
$$
(ii) For any $\beta\in \bR$,
$$
\|\Delta^{\alpha/2}u\|_{H^{\beta}_p}\leq c(\alpha,\beta)\|u\|_{H^{\beta+\alpha}_p}.
$$
\end{corollary}
\begin{remark}
                       \label{remark multiplier}
Let $\gamma, \beta \geq 0$. Then due to the well-known inequality
$$
\|u\|_{L_p}\leq \varepsilon \|u\|_{H^{\gamma}_p}+c(\varepsilon,\gamma,\beta)\|u\|_{H^{-\beta}_p},
$$
it also follows
$$
\|u\|_{H^{\gamma}_p}\leq c(\gamma,\beta, p)(\|u\|_{H^{-\beta}_p}+\|\partial^{\gamma/2}_xu\|_{L_p}).
$$
\end{remark}

For $\ell_2$-valued $\overline{\cP}$-measurable processes
$g=(g^1,g^2,\dots)$,
 we write $g\in
\bH^{\gamma}_p(T,\ell_2)$ if
\begin{equation}   \label{e:2.3}
\|g\|_{\bH^{\gamma}_p(T,\ell_2)}:= \left( \E\int^T_0\|\,
|(1-\Delta)^{\gamma/2}g(\omega, t)|_{\ell_2} \,
 \|^p_p\,dt
\right)^{1/p}<\infty.
\end{equation}
Denote $\bL_p(T) :=\bH^{0}_p (T)$ and $\bL_p(T,\ell_2)=\bH^0_p(T,\ell_2)$. Finally, we say $u_0\in U^{\gamma}_p$ if $u_0$ is $\cF_0$-measurable function $\Omega \to H^{\gamma}_p$
and
$$
\|u_0\|_{U^{\gamma}_p}:=\left( \E  \left[ \|u_0\|^p_{H^{\gamma}_p} \right]\right)^{1/p}<\infty.
$$

\begin{remark} {\rm
It is easy to check (see Remark 3.2 in \cite{Kr99} for detailed proof) that  for any
$\gamma \in (-\infty, \infty)$, $g\in \bH^{\gamma}_p(T,\ell_2)$
 and $\phi\in C^{\infty}_0(\bR^d)$ we have $\sum_{k=1}^{\infty}\int^T_0(g^k(\omega, t) ,\phi)^2 dt<\infty$ a.s., and consequently the series of stochastic integral  $\sum_{k=1}^{\infty} \int^t_0 (g^k(\omega, s),\phi)dW^k_s$ converges uniformly in $t$ in probability on $[0,T]$.
 }
\end{remark}

\begin{defn}
                                        \label{definition 3.16.1}
Write $u \in \cH^{\gamma+\alpha}_p(T)$ if $u\in
\bH^{\gamma+\alpha}_p(T), u(0)\in U^{\gamma+\alpha-\alpha/p}_p$, and
for some $f\in \bH^{\gamma}_p(T)$ and
 $ g\in \bH^{\gamma+\alpha/2}_p(T,\ell_2)$
$$
du=f dt +\sum_{k=1}^\infty g^k dW^{k}_t, \quad \hbox{for }  t\in [0,  T]
$$
 in the sense of distributions, that is, for any $\phi\in C^\infty_0 (\bR^d)$,
\begin{equation}
                \label{eqn 3.16.1}
(u(t),\phi)=(u(0),\phi)+\int^t_0(f (s),\phi)ds +
\sum_{k=1}^{\infty}\int^t_0(g^k (s),\phi)dW^{k}_s
\end{equation}
holds for all $t\leq T$ $a.s.$.
In this case we write
$$
\bD u :=f, \quad \bS^k u: =g^k, \quad \bS u :=(\bS^1 u, \dots, \bS^k u, \dots)
$$
and define the norm
\begin{equation}
                             \label{norm definition}
\|u\|_{\cH^{\gamma+\alpha}_p(T)}:=\|u\|_{\bH^{\gamma+\alpha}_p(T)}+\|\bD u\|_{\bH^{\gamma}_p(T)} +
 \|\bS u\|_{\bH^{\gamma+\alpha/2}_p(T,\ell_2)}
+\|u(0)\|_{U^{\gamma+\alpha-\alpha/p}_p}.
\end{equation}
\end{defn}

\begin{thm}
                       \label{thm banach conti}
  The space $\cH^{\gamma+\alpha}_p(T)$ is a Banach space, and
for every $ 0<t \le T$
\begin{equation}
                      \label{eqn banach1}
  \E \Big[\sup_{s\leq t}\|u(s,\cdot)\|^p_{H^{\gamma}_p} \Big] \leq
  c(p,T,\alpha)\left(\|\bD u\|^p_{\bH^{\gamma}_p(t)}+\|\bS u\|^p_{\bH^{\gamma}_p(t,\ell_2)}+\|u(0)\|^p_{U^{\gamma}_p}\right).
\end{equation}
In particular, for any $t\leq T$,
\begin{equation}
                      \label{gronwall}
\|u\|^p_{\bH^{\gamma}_p(t)}\leq c(p,T,\alpha)\int^t_0\|u\|^p_{\cH^{\gamma+\alpha}_p(s)}\,ds.
\end{equation}
\end{thm}

\begin{remark} \rm
Note that $\alpha$ is not involved in (\ref{eqn banach1}).
\end{remark}

\pf
See   Theorem 3.7 in \cite{Kr99}.  Actually in \cite{Kr99} the theorem is proved only for $\alpha=2$, but the proof works  for any $\alpha \in (0,2)$.
We will give the detailed proof of Theorem \ref{theorem banach} below, which is the  counterpart of Theorem \ref{thm banach conti} for pure-jump L\'evy processes.
\qed

\begin{remark}   \label{remark 11.22} \rm
\rm It follows from \eqref{e:2.2} that for any $\mu, \gamma\in \bR$,
the operator $(1-\Delta)^{\mu/2}:H^{\gamma}_p\to H^{\gamma-\mu}_p$
is an isometry. Indeed,
$$
\|(1-\Delta)^{\mu/2}u\|_{H^{\gamma-\mu}_p}
=\|(1-\Delta)^{(\gamma-\mu)/2}(1-\Delta)^{\mu/2}u\|_p
=\|(1-\Delta)^{\gamma/2}u\|_p=\|u\|_{H^{\gamma}_p}.
$$
The same reason shows that $(1-\Delta)^{\mu/2}:\cH^{\gamma}_p(T)\to \cH^{\gamma-\mu}_p(T)$ is an isometry.
\end{remark}


\begin{thm}
                  \label{thm deterministic}
(i) For any deterministic functions $f=f(t,x)$ and $u_0=u_0(x)$  with
$$\int^T_0
\,\|f(t, \cdot)\|^p_{H^{\gamma}_p}\,dt < \infty, \quad
\|u_0\|_{H^{\gamma+\alpha-\alpha/p}_p} < \infty,
$$
  the (deterministic) equation
$$
u_t=\Delta^{\alpha/2}u+f, \quad u(0)=u_0
$$
has a unique solution $u$
with $\int^T_0
\,\|u(t, \cdot)\|^p_{H^{\gamma+\alpha}_p}\,dt < \infty$,
 and for every $0<t \le T$
\begin{equation}\label{e:deterministic}
\int^t_0
\,\|u(s, \cdot)\|^p_{H^{\gamma+\alpha}_p} ds\leq c(p,T)\left(\int^t_0
\,\|f(s, \cdot)\|^p_{H^{\gamma}_p}\,ds
+\|u_0\|^p_{H^{\gamma+\alpha-\alpha/p}_p}\right).
\end{equation}

(ii) For any $f\in \bH^{\gamma}_p(T)$ and $u_0\in U^{\gamma+\alpha-\alpha/p}_p$, the equation
\begin{equation}
                            \label{eqn 5.29.4}
u_t=\Delta^{\alpha/2}u+f, \quad u(0)=u_0
\end{equation}
has a unique solution $u\in \bH^{\gamma+\alpha}_p(T)$ and for every $0<t \le T$
\begin{equation}
                         \label{eqn 5.20.6}
\|u\|_{\cH^{\gamma+\alpha}_p(t)}\leq c(p,T)\left(\|f\|_{\bH^{\gamma}_p(t)}
+\|u_0\|_{U^{\gamma+\alpha-\alpha/p}_p}\right).
\end{equation}

\end{thm}
\pf
(i). See,
for instance, Theorem 2.1 in \cite{M1}.

(ii). This result is also known. See, for instance, Lemma 3.2 and Lemma 3.4 of  \cite{CL}.   Actually since  equation (\ref{eqn 5.29.4}) is deterministic for each fixed $\omega$, the claim of (ii) can be  obtained from (i).  Indeed, the uniqueness and
estimate (\ref{eqn 5.20.6}) are obvious by (i). For the existence of solution,
assume that $u_0$ and $f$ are sufficiently smooth in $x$, then using Fourier transform one can easily check that
$$
u(t):=T_tu_0 +\int^t_0T_{t-s}f ds
$$
solves (\ref{eqn 5.29.4}) and is in $\cH^{\gamma+\alpha}_p(T)$.
 For general $u_0$ and $f$ it is enough to use a standard approximation argument (see, for instance, the proof Theorem \ref{thm linear}).

\qed

Now we give our assumption on $a(\omega,t)$.
\begin{assumption}
                  \label{ass 1}
The process $a(\omega,t)$ is predictable and there is a constant $\delta>0$ so that
$$
\delta<a(\omega,t)<\delta^{-1}, \quad \quad \forall \omega,t.
$$
\end{assumption}

Now we present an $L_p$-theory for linear stochastic parabolic equations with random fractional Laplacian.

\begin{thm}
                        \label{thm linear}
Let $p\in [2,\infty)$ and $\gamma\in \bR$. For any $f\in \bH^{\gamma}_p(T)$, $g\in \bH^{\gamma+\alpha/2}_p(T,\ell_2)$ and $u_0\in U^{\gamma+\alpha-\alpha/p}_p$, the linear equation
\begin{equation}
                                \label{eqn main}
du=\left(a(\omega,t)\Delta^{\alpha/2}u+f\right)\,dt
+\sum_{k=1}^{\infty}g^kdW^k_t, \quad u(0)=u_0,
\end{equation}
admits a unique solution $u$ in $\cH^{\gamma+\alpha}_p(T)$, and for this solution
\begin{equation}
                     \label{eqn linear}
\|u\|_{\cH^{\gamma+\alpha}_p(T)}\leq c(p,T,\delta)\left(\|f\|_{\bH^{\gamma}_p(T)}
+\|g\|_{\bH^{\gamma+\alpha/2}_p(T,\ell_2)}+\|u_0\|_{U^{\gamma+\alpha-\alpha/p}_p}\right).
\end{equation}
\end{thm}

\begin{remark}
                  \label{remark kijung}
(i) Recall that the unique solution $u\in \cH^{\gamma+\alpha}_p(T)$ is understood in the sense of distributions as in Definition \ref{definition 3.16.1}, that is, for any  $\phi\in C^\infty_0 (\bR^d)$,
$$
(u(t),\phi)=(u(0),\phi)+\int^t_0 \left(a(\omega,s)(u,\Delta^{\alpha/2}\phi)+(f(s),\phi)\right) ds +
\sum_{k=1}^{\infty}\int^t_0(g^k (s),\phi)dW^{k}_s
$$
holds for all $t\leq T$ $a.s.$.

(ii)   A version of Theorem \ref{thm linear} is proved in \cite{CL} under stronger conditions on $g$ and the processes. Precisely in \cite{CL} it is assumed that  $a(\omega,t)=1$, $g\in \bH^{\gamma+\alpha/2+\varepsilon}_p(T)$, $\varepsilon>0$,
and there are only finitely many Wiener processes in  equation (\ref{eqn main}).
\end{remark}

\pf  {\bf Step 1}. Owing to Remark \ref{remark 11.22},  we only need to
show that the theorem holds for a particular $\gamma=\gamma_0$.
Indeed, suppose that the theorem holds when $\gamma=\gamma_0$. Then
it is enough to notice that $u\in \cH^{\gamma+\alpha}_p(T)$ is a
solution of the equation if and only if
$\bar{u}:=(1-\Delta)^{(\gamma-\gamma_0)/2}u\in \cH^{\gamma_0+\alpha}_p$ is
a solution of the equation with
$$
\bar{f}:=(1-\Delta)^{\frac{(\gamma-\gamma_0)}{2}}f,\quad
\bar{g}:=(1-\Delta)^{\frac{(\gamma-\gamma_0)}{2}}g, \quad
 \bar{u}_0:=(1-\Delta)^{\frac{(\gamma-\gamma_0)}{2}}u_0,
 $$ in place of
$f,g$ and $u_0$, respectively. Furthermore,
\begin{eqnarray*}
\|u\|_{\bH^{\gamma+\alpha}_p(T)}=\|\bar{u}\|_{\bH^{\gamma_0+\alpha}_p(T)}
&\leq& c\left(\|\bar{f}\|_{\bH^{\gamma_0}_p(T)}
+\|\bar{g}\|_{\bH^{\gamma_0+\alpha/2+\varepsilon_0}_p(T, \ell_2)}
+\|\bar{u}_0\|_{U^{\gamma_0+\alpha-\alpha/p}_p} \right) \\
&=&c \left(
\|f\|_{\bH^{\gamma}_p(T)}+\|g\|_{\bH^{\gamma+\alpha/2+\varepsilon_0}_p(T,
\ell_2)}+\|u_0\|_{U^{\gamma+\alpha-\alpha/p}_p} \right).
\end{eqnarray*}

{\bf Step 2}. Next we  assume $a(\omega,t)=1$ and  prove the theorem for the equation:
\begin{equation}
                        \label{heat eqn}
  du=\Delta^{\alpha/2} u dt + \sum_{k=1}^{\infty} g^k dW^k_t, \quad u(0)=0.
  \end{equation}
Remember that we may assume   $\gamma=-\alpha/2$.  Since the uniqueness of (\ref{heat eqn}) follows from results for the deterministic equations (Theorem \ref{thm deterministic}), we only need to show that there exists a solution $u\in \bH^{\alpha/2}_p(T)$
of \eqref{heat eqn}
 and $u$ satisfies estimate (\ref{eqn linear}) with $f=u_0=0$ and $\gamma=-\alpha/2$.

 For a moment,  assume
$N_0>0$ is a fixed non-random constant, $g^k=0$ for all $k> N_0$ and
\begin{equation}
                       \label{eqn 03.17.11}
g^k(t,x)=\sum_{i=0}^{m_k}I_{(\tau^k_i,\tau^k_{i+1}]}(t)g^{k_i}(x) \qquad \text{for } k \le N_0,
\end{equation}
where $\tau^k_i$ are bounded stopping times and $g^{k_i}(x)\in
 C^{\infty}_0 (\bR^d)$. Define
$$
v(t,x):=\sum_{k=1}^{N_0} \int^t_0
 g^k(s, x) dW^k_s
=\sum_{k=1}^{N_0}\sum_{i=1}^{m_k} g^{k_i}(x)(W^k_{t\wedge
\tau^k_{i+1}}-W^k_{t\wedge \tau^k_{i}})
$$
and
\begin{equation}
                    \label{eqn 03-22-1}
u(t,x):=v(t,x)+\int^t_0 \Delta^{\alpha/2}T_{t-s} v (s,x)\,
ds=v(t,x)+\int^t_0 T_{t-s}\Delta^{\alpha/2} v (s,x)\, ds.
\end{equation}
 Using Fourier transform one can easily   show (See, for instance, \cite{CL}) that if  functions $h_1=h_1(t,x)$ and $h_2=h_2(x)$ are  sufficiently smooth in $x$ then
$$
w_1(t,x):=\int^t_0T_{t-s}h_1(s,x) ds, \quad w_2(t,x)=T_t h_2 (x)
$$
solve
$$
dw_1=(\Delta^{\alpha/2}w_1 + h_1)\,dt, \quad w_1(0)=0,
$$
$$
dw_2=\Delta^{\alpha/2}w_2\,dt, \quad w_2(0)=h_2.
$$
Therefore we have $d(u-v)=(\Delta^{\alpha/2}(u-v)+\Delta^{\alpha/2}
v)dt=\Delta^{\alpha/2} u dt$, and
$$
du=\Delta^{\alpha/2} u dt +dv=\Delta^{\alpha/2} u dt + \sum_{k=1}^{N_0} g^kdW^k_t.
$$
Also by (\ref{eqn 03-22-1}) and stochastic Fubini theorem (\cite[Theorem 64]{P}), almost
surely,
\begin{eqnarray}\label{e:fdsa}
u(t,x)&=&v(t,x)+\sum_{k=1}^{N_0}\int^t_0 \int^s_0
\Delta^{\alpha/2}T_{t-s}
g^k (r, x)dW^k_r ds \nonumber\\
&=&v(t,x)-\sum_{k=1}^{N_0}\int^t_0\int^t_r\frac{\partial}{\partial
s}T_{t-s}g^k (r, x) ds dW^k_r\nonumber\\
 &=&\sum_{k=1}^{N_0} \int^t_0 T_{t-s}g^k (s,x) dW^k_s.
\end{eqnarray}
Hence,
$$
\partial^{\alpha/2}_xu(t,x)=\sum_{k=1}^{N_0}\int^t_0
\partial^{\alpha/2}_xT_{t-s}g^k (s,\cdot) (x)dW^k_s,
$$
and by
 Burkholder-Davis-Gundy's
inequality, we have
$$
\E \left[  \big|\partial^{\alpha/2}_xu(t,x)\big|^p \right] \leq c(p)\E\left[
\left(\int^t_0 \sum_{k=1}^{N_0}
|\partial^{\alpha/2}_xT_{t-s}g^k(s,\cdot) (x)|^2ds\right)^{p/2}\right].
$$
Now we use a parabolic version of Littlewood-Paley inequality for fractional Laplacian (Theorem 2.3 in \cite{Ildoo})
\begin{equation}
                    \label{eqn ildoo}
\int_{\bR^d}\int^T_0\left
[\int^t_0|\partial^{\alpha/2}_xT_{t-s}g(s,\cdot)(x)|^2_{\ell_2}
ds\right]^{p/2}dtdx\leq c(\alpha,p)\int_{\bR^d}\int^T_0
|g(t,x)|^p_{\ell_2}\,dtdx
\end{equation}
and get
\begin{equation}
                               \label{e:12.2.3}
\E \left[ \int^T_0\|\partial^{\alpha/2}_xu(t, \cdot)\|^p_{p} \, dt \right]\leq c(p)\E
\left[ \int^T_0\||g(t, \cdot)|_{\ell_2}\|^p_{p} \, dt\right].
\end{equation}
Similarly,
 \eqref{e:fdsa} and Burkholder-Davis-Gundy's
inequality yield
\begin{equation}
                                 \label{eqn 8.6.5}
\E \left[ |u(t,x)|^p \right] \leq c(p)\, \E\left[ \left(
\int^t_0\sum_{k=1}^{N_0}|T_{t-s}g^k(s,x)|^2ds\right)^{p/2} \right].
 \end{equation}
  Since $(\sum_{k=1}^{N_0}|a_n|^2)^{p/2}\leq
 c(N_0,p)\sum_{k=1}^{N_0}|a_n|^p$ and $\|T_tf\|_p\leq
 c\|f\|_p$,
 we see that for every $t>0$
\begin{align*}
 \int_{\bR^d}\left( \int^t_0\sum_{k=1}^{N_0}|T_{t-s}g^k(s,x)|^2ds\right)^{p/2} dx
 &\le  t^{p/2-1} \int_{\bR^d}\int^t_0  \left(\sum_{k=1}^{N_0}|T_{t-s}g^k(s,x)|^2           \right)^{p/2}  dt dx \\
 &\le c(T,N_0,p)\int^t_0\int_{\bR^d}\sum_{k=1}^{N_0}|g^k(t,x)|^p dxdt.
 \end{align*}
 Consequently,
 \begin{align}\label{e:hjkrty}
&\E \int_0^T\int_{\bR^d} |u(t,x)|^p  dxdt   \le c(T,N_0,p)  \E \int^T_0  \int_{\bR^d} |g(t,x)|^p_{\ell_2}  dx     ds.
  \end{align}
 Thus we proved $\partial^{\alpha/2}_x u, u\in \bL_{p}(T)$, and  hence by Corollary \ref{cor multiplier} we have  $u\in \cH^{\alpha/2}_p(T)$.  Note that by Corollary \ref{cor multiplier}(ii)
 $$
 \|\Delta^{\alpha/2}u\|_{H^{-\alpha/2}_p}=\|\Delta^{\alpha/4}(\partial^{\alpha/2}_xu)\|_{H^{-\alpha/2}_p}\leq c \|\partial^{\alpha/2}_xu\|_{L_p}.
 $$
  By definition (\ref{norm definition}) and  Remark \ref{remark multiplier}, for any $t\leq T$,
 \begin{eqnarray}
 \|u\|^p_{\cH^{\alpha/2}_p(t)}&\leq& c(p)\left(\|u\|^p_{\bH^{\alpha/2}_p(t)}+\|\Delta^{\alpha/2}u\|^p_{\bH^{-\alpha/2}_p(t)}+\|g\|^p_{\bL_p(t,\ell_2)}\right)\nonumber\\
 &\leq& c\left(\|u\|^p_{\bH^{-\alpha/2}_p(t)}+\|\partial^{\alpha/2}_xu\|^p_{\bL_p(t)}+\|g\|^p_{\bL_p(t,\ell_2)}\right). \label{e1}
 \end{eqnarray}
  Combining this with  (\ref{e:12.2.3}) and (\ref{gronwall})
 we have that for every $0<t \le T$
  \begin{eqnarray}
\|u\|^p_{\cH^{\alpha/2}_p(t)}&\leq& c(p,T,\alpha)\left(\|u\|^p_{\bH^{-\alpha/2}_p(t)} +
\|g\|^p_{\bL_p(t,\ell_2)} \right)\nonumber\\
&\leq& c(p,T,\alpha) \int_0^t\|u\|^p_{\cH^{\alpha/2}_p(s)} ds +  c(p,T,\alpha)
\|g\|^p_{\bL_p(T,\ell_2)}. \label{e2}
 \end{eqnarray}
 Finally,  Gronwall  leads to
  (\ref{eqn linear}).

 Now we drop the additional assumptions on $g$ by using the following standard  approximation argument:
    By Theorem 3.10 in \cite{Kr99}, for $g\in \bL_p(T,\ell_2)$ we can take a sequence
  $g_n \in \bL_p(T,\ell_2)$ so that $g_n\to g$ in $\bL_p(T,\ell_2)$ and
  each $g_n=(g^1_n,g^2_n,\cdots)$ satisfies above assumed assumptions, that is,
  $g^k_n=0$ for all large $k$ and each $g^k_n$ is of type (\ref{eqn 03.17.11}).
  By the above result, the equation
 $$
 du_n=\Delta^{\alpha/2}u_ndt+ \sum_{k=1}^\infty g^k_n dW^k_t, \quad
 u_n(0)=0
 $$
 has a unique solution $u_n$. It also follows that    $u_n-u_m$ is the unique solution of
 $$
d(u_n-u_m)=\Delta^{\alpha/2}(u_n-u_m)dt+ (g^k_n-g^k_m) dW^k_t, \quad  
 (u_n-u_m)(0)=0
$$
and, by the previous argument
$$
\|u_n-u_m\|_{\cH^{\alpha/2}_p(T)}\leq c(p,T)\|g_n-g_m\|_{\bL_p(T,\ell_2)}.
$$
Consequently, there is $u\in \cH^{\alpha/2}_p(T)$ so that $u_n \to u$ in $\cH^{\alpha/2}_p(T)$.
We only need to  prove $u$ is a solution of (\ref{heat eqn}).  Equivalently, we need to prove that for any
$\phi\in C^{\infty}_0 (\bR^d)$,  the equality
\begin{equation}
                            \label{eqn 11111}
(u(t, \cdot),\phi)=\int^t_0(\Delta^{\alpha/2}u(s, \cdot),\phi)ds+ \sum_{k=1}^{\infty} \int^t_0 (g^k(s, \cdot),\phi)dW^k_s.
\end{equation}
holds for for all $t\leq T$ (a.s.), or equivalently
\begin{equation}
                        \label{eqn 6.1.1}
((1-\Delta)^{-\alpha/2}u(t, \cdot), (1-\Delta)^{\alpha/2}\phi)=\int^t_0(\Delta^{\alpha/4}u(s, \cdot),\Delta^{\alpha/4}\phi)ds+ \sum_{k=1}^{\infty} \int^t_0 (g^k(s, \cdot),\phi)dW^k_s.
\end{equation}
By  (\ref{eqn banach1}),
$$
\lim_{n \to \infty} \E \left[\sup_{t\leq T}\|(1-\Delta)^{-\alpha/2}(u_n(t, \cdot)-u(t, \cdot))\|^p_{L_p} \right]=0, \quad \text{a.s.}
$$
which implies that one can take a subsequence $n_j$ so that
 $(1-\Delta)^{-\alpha/2}u_{n_j}\to (1-\Delta)^{-\alpha/2}u$ in $L_p(\bR^d)$ uniformly on $[0,T]$ (a.s) and consequently $t \to ((1-\Delta)^{-\alpha/2}u(t, \cdot ),(1-\Delta)^{\alpha/2}\phi)$ is continuous  on $[0,T]$.
By taking the limit from
$$
((1-\Delta)^{-\alpha/2}u_{n_j}(t, \cdot ),(1-\Delta)^{\alpha/2}\phi)=\int^t_0(\Delta^{\alpha/4}u_{n_j}(s, \cdot ),\Delta^{\alpha/4}\phi)ds+ \sum_{k=1}^{\infty}\int^t_0 (g^k_{n_j}(s, \cdot ),\phi)dW^k_s
$$
and remembering that both sides of (\ref{eqn 6.1.1}) are continuous in $t$, one easily get that equality
(\ref{eqn 6.1.1}) holds for all $t\leq T$ (a.s.).

{\bf{Step 3}}. Next we prove the theorem for the equation
\begin{equation}
                        \label{heat eqn2}
  du=(\Delta^{\alpha/2} u +f) dt +  g^k dW^k_t, \quad u(0)=u_0.
  \end{equation}
  Again we may assume $\gamma=-\alpha/2$, and due to Theorem \ref{thm deterministic} we only need to show that there exists a solution
  $u$ and it satisfies  estimate (\ref{eqn linear}). By Theorem \ref{thm deterministic}, the equation
  $$
  dv=(\Delta^{\alpha/2}v +f)dt, \quad v(0)=u_0
  $$
  has a solution $v\in \cH^{\alpha/2}_p(T)$ and
  $$
  \|v\|_{\cH^{\alpha/2}_p(T)}\leq c(p,T)\left(\|f\|_{\bH^{-\alpha/2}_p(T)}+\|u_0\|_{U^{\alpha/2-\alpha/p}_p}\right).
  $$
  Also by the result of {\bf Step 2}, the equation
  $$
  dw=\Delta^{\alpha/2} wdt +  \sum_{k=1}^\infty g^k dW^k_t, \quad w(0)=0
  $$
  has a unique solution and
  $$
  \|w\|_{\cH^{\alpha/2}_p(T)}\leq c(p,T)\|g\|_{\bL_p(T,\ell_2)}.
  $$
  Now it is enough to take $u=v+w$.

  {\bf{Step 4}} (A priori estimate). We prove the a priori estimate (\ref{eqn linear}) holds given that a solution $u\in \cH^{\gamma+\alpha}_p(T)$ of  the following equation already exists :
  $$
du=\left(a(\omega,t)\Delta^{\alpha/2}u+f\right)\,dt
+\sum_{k=1}^{\infty}g^kdW^k_t, \quad u(0)=u_0.
$$
  This time we prove (\ref{eqn linear}) only for  $\gamma=0$. This is enough due to the reason given  in Step 1.  By  Step 3, the equation
  $$
  dv=(\Delta^{\alpha/2} v +f) dt +  \sum_{k=1}^\infty g^k dW^k_t, \quad v(0)=u_0.
  $$
  has a solution $v\in \cH^{\alpha}_p(T)$ and
  $$
   \|v\|_{\cH^{\alpha}_p(T)}\leq c(p,T)\left(\|f\|_{\bL_p(T)}+\|g\|_{H^{\alpha/2}_p(T,\ell_2)}+\|u_0\|_{U^{\alpha-\alpha/p}_p}\right).
  $$
  Note that $\bar{u}:=u-v$ satisfies
  $$
  d\bar{u}=(a(\omega,t)\Delta^{\alpha/2}\bar{u} +\bar{f})dt, \quad \bar{u}(0)=0,
  $$
  where $\bar{f}:=(a(\omega,t)-1)\Delta^{\alpha/2}v$, and
  $$
  \|\bar{f}\|_{\bL_p(T)}\leq c\|\Delta^{\alpha/2}v\|_{\bL_p(T)}\leq c \left(\|f\|_{\bL_p(T)}+\|g\|_{H^{\alpha/2}_p(T,\ell_2)}+\|u_0\|_{U^{\alpha-\alpha/p}_p}\right).
  $$
  Since $u=v+\bar{u}$,  $\|u\|_{\cH^{\alpha}_p(T)}\leq \|v\|_{\cH^{\alpha}_p(T)}+\|\bar{u}\|_{\cH^{\alpha}_p(T)}$ and $\|\bar{u}\|_{\cH^{\alpha}_p(T)}\leq c\|\bar{u}\|_{\bH^{\alpha}_p(T)}+c\|\bar{f}\|_{\bL_p(T)}$,
  to prove (\ref{eqn linear}) we only need to show that  for each   $\omega\in \Omega$,
  \begin{equation}
                    \label{eqn 3.17.7}
   \int^T_0 \|\bar{u}(t, \cdot)\|^p_{H^{\alpha}_p}dt \leq c(p,T,\alpha,\delta)\int^T_0 \|\bar{f}(t, \cdot)\|^p_{L_p}dt.
  \end{equation}
 For fixed $\omega$, define a non-random
  functions
  \begin{equation}
                          \label{relation}
  \tilde{u}(t,x)=\bar{u}(\omega,\xi (\omega, t),x)
   \quad \tilde{f}(t,x)=a(\omega, t)^{-1} \bar{f}(\omega,\xi (\omega, t),x).
    \end{equation}
    where
 $ \xi (\omega, t):=  \int^t_0 \frac{ds}{a(\omega, s)}$.
    Then
    clearly
    $\tilde{u}$ satisfies
  $$
  \tilde{u}_t=\Delta^{\alpha/2}\tilde{u} + \tilde{f}, \quad \tilde{u}(0)=0.
  $$
Let $\tilde  T(\omega, T)$ be such that $T= \int^{ \tilde  T(\omega, T)}_0 \frac{ds}{a(\omega, s)}$.
 Since  $   \delta  T <\tilde  T(\omega, T)< \delta T$, applying \eqref{e:deterministic}, we get
$$
   \int^{\tilde  T(\omega, T)}_0 \|\tilde{u}(t, \cdot)\|^p_{H^{\alpha}_p}dt \leq c(p,T,\delta,\alpha)\int^{\tilde  T(\omega, T)}_0 \|\tilde{f}(t, \cdot)\|^p_{L_p}dt.
$$
This and relations in (\ref{relation}) easily lead to (\ref{eqn 3.17.7}).

  {\bf{Step 5}} (Method of continuity). The solvability of equation (\ref{heat eqn2}),  the a priori estimate (\ref{eqn linear}) and the method of continuity obviously finish the proof of the theorem. But  below we  show how the method of continuity works  only for  reader's convenience.

  For $\lambda\in [0,1]$, denote $a_{\lambda}(\omega,t)=(1-\lambda)+\lambda a(\omega,t)$. Then obviously $a_{\lambda}$ is predictable and $a_{\lambda}\in (\delta,\delta^{-1})$. It follows from Step 4 that if $u\in \cH^{\gamma+\alpha}_p(T)$  is a solution of  the equation
   \begin{equation}
                          \label{eqn 03.18.8}
  du=\left(a_{\lambda}(\omega,t)\Delta^{\alpha/2}u+f\right)\,dt
+\sum_{k=1}^{\infty}g^kdW^k_t, \quad u(0)=u_0,
\end{equation}
  then  the estimate (\ref{eqn linear}) holds with the same constant $c=c(p,T,\delta)$. Now let $J$ be the collection of $\lambda \in [0,1]$ so that for any $f\in \bH^{\gamma}_p(T), g\in \bH^{
  \gamma
  +\alpha/2}_p(T,\ell_2)$ and $u_0\in U^{\gamma+\alpha/2-\alpha/p}_p$,
  equation (\ref{eqn 03.18.8}) has a solution.  By Step 3, $0\in J$. Note that to finish the proof of the theorem we only need to show $1\in J$.   Now let $\lambda_0\in J$. Obviously $u\in \cH^{\gamma+\alpha}_p(T)$ is a solution of (\ref{eqn 03.18.8}) if and only if
  \begin{equation}
                          \label{eqn 03.18.9}
  du=\left(a_{\lambda_0}(\omega,t)\Delta^{\alpha/2}u+[(a_{\lambda}(\omega,t)-a_{\lambda_0}(\omega,t))\Delta^{\alpha/2}u+f]\right)\,dt
+\sum_{k=1}^{\infty}g^kdW^k_t, \quad u(0)=u_0.
\end{equation}
  Now fix $u^1 \in  \cH^{\gamma+\alpha}_p(T)$ with initial date $u_0$ (for instance take the solution of (\ref{heat eqn2})), and define $u^2, u^3,\cdots$ so that $u^{n+1}\in\cH^{\gamma+\alpha}_p(T)$ is the solution of
  \begin{equation}
                          \label{eqn 03.18.10}
  du^{n+1}=\left(a_{\lambda_0}(\omega,t)\Delta^{\alpha/2}u^{n+1}+[(a_{\lambda}(\omega,t)-a_{\lambda_0}(\omega,t))\Delta^{\alpha/2}u^n+f]\right)dt\,dt
+ g^kdW^k_t, \quad u(0)=u_0.
\end{equation}
Then $v^{n+1}:=u^{n+1}-u^n$ satisfies
$$
 dv^{n+1}=\left(a_{\lambda_0}(\omega,t)\Delta^{\alpha/2}v^{n+1}+(a_{\lambda}(\omega,t)-a_{\lambda_0}(\omega,t))\Delta^{\alpha/2}v^n\right)\,dt
$$
By the a priori estimate (\ref{eqn linear}),
\begin{eqnarray*}
\|v^{n+1}\|_{\cH^{\gamma+\alpha}_p(T)} &\leq& c\|(a_{\lambda}-a_{\lambda_0})\Delta^{\alpha/2}v^n\|_{\bH^{\gamma}_p(T)}\\
&\leq& N(p,T,\delta) |\lambda-\lambda_0|\|v^{n}\|_{\cH^{\gamma+\alpha}_p(T)}.
\end{eqnarray*}
Thus if $|\lambda-\lambda_0|<1/(2N(p,T,\delta))$, the map  which send $u^n$ to  $u^{n+1}$ is a contraction in
 $\cH^{\gamma+\alpha}_p(T)$, and  has a unique fixed point $u$. Thus $u$ satisfies
 (\ref{eqn 03.18.8})--(\ref{eqn 03.18.9}).  Since the above constant $N$ is independent of $\lambda$, it follows that $J=[0, 1]$ and the theorem is proved.
\qed

Finally we consider the nonlinear equation
\begin{equation}
                                \label{eqn nonlinear}
du=\left(a(\omega,t)\Delta^{\alpha/2}u+f(u)\right)\,dt
+\sum_{k=1}^{\infty}g^k(u)dW^k_t, \quad u(0)=u_0,
\end{equation}
where $f(u)=f(\omega,t,x,u)$ and $g^k(u)=g^k(\omega,t,x,u)$.

\begin{assumption}
                 \label{ass nonlinear}
           Assume       $f(0) \in \bH^{\gamma}_p(T)$ and $g(0) \in \bH^{\gamma+\alpha/2}_p(T, \ell_2)$.
Moreover,   for any $\varepsilon>0$, there exists a constant $K_{\varepsilon}$ so that
  for any $u=u(x),v=v(x) \in H^{\gamma+\alpha}_p$ and
  $\omega,t$,
  we have
  \begin{align}
                          \label{eqn 03.18.91}
 & \|f(t, \cdot, u(\cdot))-f(t, \cdot, v(\cdot))\|_{H^{\gamma}_p} +\|g(t, \cdot, u(\cdot))-g(t, \cdot, v(\cdot))\|_{H^{\gamma+\alpha/2}_p( \ell_2)}
 \nonumber\\
 \leq &\varepsilon \|u-v\|_{H^{\gamma+\alpha}_p}+K(\varepsilon)\|u-v\|_{H^{\gamma}_p}.
  \end{align}
\end{assumption}

To give an example of $f(u)$ and $g(u)$ satisfying Assumption \ref{ass nonlinear}, we introduce the space of point-wise multipliers in $H^{\gamma}_p$.  For each $r\geq 0$,  define
 \begin{equation}\label{e:2.11}
 B^{r}=\begin{cases} B(\bR^d)  \qquad & \hbox{if }  r=0,\\
 C^{r-1,1}(\bR^d) &\hbox{if }  r=1,2, \cdots ,\\
C^{r}(\bR^d) &  \hbox{otherwise},
 \end{cases}
\end{equation}
where $B(\bR^d)$ is the space of
bounded Borel measurable functions
on $\bR^d$, $C^{r-1,1}(\bR^d)$ is the space of $r-1$ times
continuously differentiable functions whose $(r-1)$st order derivatives  are Lipschitz continuous, and
$C^{r}(\bR^d)$ is the usual H\"older space. Also we use the space $B^{r}$ for $\ell_2$-valued functions.
For instance,
 if $g=(g^1,g^2,...)$,  then $|g|_{B^0}=\sup_x |g(x)|_{\ell_2}$ and
$$
|g|_{C^{n-1,1}}=\sum_{|\alpha|\leq n-1} |D^{\alpha}g
|_{B^0}+\sum_{|\alpha|=n-1}\sup_{x\neq y}\frac{
|D^{\alpha}g(x)-D^{\alpha}g(y)|_{\ell_2}}{|x-y|}.
$$
 Fix $\kappa_0=\kappa_0(\gamma)\geq 0$ so that $\kappa_0>0$ if $\gamma$ is not integer.
It is  known  (see, for instance, Lemma 5.2 in \cite{Kr99})  that for any $a\in B^{|\gamma|+\kappa_0}$ and $h\in H^{\gamma}_p$,
\begin{equation}
                                             \label{multiplier}
\|ah\|_{H^{\gamma}_p}\leq c(\gamma,\kappa_0)|a|_{B^{|\gamma|+\kappa_0}}|h|_{H^{\gamma}_p}
\end{equation}
and the same  inequality holds for    $\ell_2$-valued functions $a$.

\begin{example}
                             \label{ex 1}
 Fix $\kappa_0=\kappa_0(\gamma)\geq 0$ so that $\kappa_0>0$ if $\gamma$ is not integer.
Consider
$$f(u)=b(\omega,t,x)\Delta^{\beta_1/2}u+ \sum_{i=1}^d c^i(\omega,t,x)u_{x^i}I_{\alpha>1}+d(\omega,t,x)u+f_0,
$$
$$
 g^k(u)=\sigma^{k}(\omega,t,x)
 \Delta^{
 \beta_2/2}u+v^k
(\omega,t,x)u+g^k_0,
 $$
 where $\beta_1<\alpha$, $\beta_2<\alpha/2$,  $f_0 \in \bH^{\gamma}_p(T)$ and  $g_0\in \bH^{\gamma+\alpha/2}_p(T,\ell_2)$.  Assume for each $\omega,t$,
 $$
 |b|_{B^{|\gamma|+\kappa_0}}+\sum_{i=1}^d|c^i|_{B^{|\gamma|+\kappa_0}}+|d|_{B^{|\gamma|+\kappa_0}}
 +|\sigma|_{B^{|\gamma|+\alpha/2+\kappa_0}}+|\nu|_{B^{|\gamma|+\alpha/2+\kappa_0}}\leq K.
 $$
Then by (\ref{multiplier}), for each $t$
\begin{eqnarray*}
&&\|f(t, \cdot, u(\cdot))-f(t, \cdot, v(\cdot))\|_{H^{\gamma}_p} +\|g(t, \cdot, u(\cdot))-g(t, \cdot, v(\cdot))\|_{H^{\gamma+\alpha/2}_p(\ell_2)}\\
&\leq& c\left(\|\Delta^{\beta_1/2}(u-v)\|_{H^{\gamma}_p}+I_{\alpha>1}\|D(u-v)\|_{H^{\gamma}_p}
+\|u-v\|_{H^{\gamma+ \alpha/2}_p}
+\|\Delta^{\beta_2/2}(u-v)\|_{H^{\gamma+\alpha/2}_p}\right).
\end{eqnarray*}
Since for any $\alpha_1<\alpha$ and $\varepsilon>0$, by interpolation theory,
$$\|u\|_{H^{\gamma+\alpha_1}_p}\leq   c(\alpha, \alpha_1) \|u\|_{H^{\gamma+\alpha}_p}^{\alpha_1/\alpha}\|u\|_{H^{\gamma}_p}^{1-\alpha_1/\alpha}
\leq \varepsilon \|u\|_{H^{\gamma+\alpha}_p}+c(\varepsilon,\alpha_1, \alpha)\|u\|_{H^{\gamma}_p},$$
one easily gets (\ref{eqn 03.18.91}).
\end{example}

Here is the main result of this section.

\begin{thm} \label{e:twcase}
Suppose Assumptions \ref{ass 1} and \ref{ass nonlinear} hold.
Then equation (\ref{eqn nonlinear}) has a unique solution $u\in \cH^{\gamma+\alpha}_p(T)$, and for this solution
$$
  \|u\|_{\cH
  ^{\gamma+\alpha}
  _p(T)}\leq c\left(\|f(0)\|_{\bH^{\gamma}_p(T)}
  +\|g(0)\|_{\bH^{\gamma+\alpha/2}_p(T,\ell_2)}
  +\|u_0\|_{U^{\gamma+\alpha/2-\alpha/p}_p}\right),
 $$
 where $c=c(p,T,\delta)$.
 \end{thm}
\pf
Our proof is virtually identical to the that of Theorem 6.4 in \cite{Kr99}, where the theorem is proved when $\alpha=2$. The only difference is that one has to use Theorem \ref{thm linear} in this article, in place the corresponding result in \cite{Kr99}.
We skip the proof here since we will give the proof for more general case in next section.
\qed

By  Sobolev embedding theorem, we immediately get the following
\begin{corollary}\label{c:reg1}
Suppose Assumptions \ref{ass 1} and\ref{ass nonlinear} hold.
If $\gamma + \alpha > d/p$, then the unique solution $u\in \cH^{\gamma+\alpha}_p(T)$ of
equation (\ref{eqn nonlinear})
 is $C^{\gamma + \alpha-d/p}$-valued process on $[0, T] \times \Omega$ a.s..
\end{corollary}

\mysection{General case}
                                    \label{section 3}

Let $Z^{1}_{t},Z^{2}_{t},...$
 be independent $m$-dimensional L\'evy processes relative to $\{\cF_{t},t\geq0\}$.
For $t\geq 0$ and Borel set $A\in \cB(\bR^m\setminus \{0\})$,  define
$$
N_k(t,A) :=\# \left\{0\leq s\leq t; \, Z^k_s-Z^k_{s-} \in A
\right\}, \quad \wt {N}_k(t,A):=N_k(t,A)-t\nu_k(A)
$$
where $\nu_k(A):=\E [N_k(1,A)]$ is the L\'e{}vy measure of $Z^k$. By
L\'e{}vy-It\^o decomposition,  there exist  a vector $\alpha^k$,  a non-negative definite matrix
$\beta^k$ and $m$-dimensional  Wiener process $B^k$ so that
\begin{equation}
                                           \label{eqn 2.27}
Z^k(t)=\alpha^kt +\beta^k B^k_t+\int_{|z|<1}z \wt
{N}_k(t,dz)+\int_{|z|\geq 1} z N_k(t, dz).
\end{equation}
For any $q,k=1,2,\cdots$, denote
$$
\wh {c}_{k,q}:=\left(\int_{\bR^m} |z|^q \nu_k(dz)\right)^{1/q}.
$$
Now we fix $p\in [2,\infty)$ and denote $\wh{c}_k:= \left(\wh {c}_{k,2}\vee \wh
{c}_{k,p}\right)$. In this section we assume

\begin{equation}
             \label{main assump}
\wh {c}:=\sup_{k \ge 1} \wh{c}_k <\infty.
\end{equation}
((\ref{main assump}) will be weaken in section \ref{section 4}).
Then for any    $2< q<p$, by H\"older's inequality,
$$ \wh c_{k, q} \leq \left(\int_{\bR^m} |z|^2 \nu_k(dz)\right)^{(p-q)/(q(p-2))}
\left(\int_{\bR^m} |z|^p \nu_k(dz)\right)^{(q-2)/(q(p-2))} \leq \wh
c_k.
$$
By (\ref{main assump}),   $\int_{|z|\geq 1} |z| \nu_k(dz)\leq
\int_{|z|\geq 1} |z|^2 \nu_k(dz)<\infty$, and
$$
\int_{|z|\geq 1} z N_k(t, dz)=\int_{|z|\geq 1}z \wt
{N}_k(t,dz)+t\int_{|z|\geq 1}z \nu_k(dz).
$$
Thus by absorbing $\wt \alpha_k:=\int_{|z|\geq 1} z \nu_k(dz)$ into
$\alpha_k$ we can rewrite (\ref{eqn 2.27}) as
$$
Z^k_t=\tilde{\alpha}_k t +\beta_k B^k_t+\int_{\bR^m}z \wt
{N}_k(t,dz).
$$
We first consider the following linear equation:
\begin{equation}
                                \label{eqn main levy_old}
du=\left(a(\omega,t)\Delta^{\alpha/2}u+f\right)dt\,
+\sum_{k=1}^{\infty}g^k \cdot dZ^k_t, \quad u(0)=u_0.
\end{equation}
Relocation of  the term $\sum_{k=1}^{\infty} g^k \cdot \tilde{\alpha}_k dt$  into the deterministic part of
(\ref{eqn main levy_old}) allow us to  assume $\tilde{\alpha}_k=(0, \dots, 0)$. Moreover, since $B^{k,j}$'s are independent 1-dimensional Wiener processes where
$B^{k}=(B^{k,1}, \dots B^{k,m})$,
\eqref{eqn main levy_old} can be written as
\begin{equation}
                                \label{eqn main levy}
du=\left(a(\omega,t)\Delta^{\alpha/2}u+f\right)\,dt
+
\sum_{i=1}^{\infty}h^k dW^k_t +
\sum_{k=1}^{\infty}\sum_{j=1}^{m}g^{k,j}  d Y^{k,j}_t, \quad u(0)=u_0,
\end{equation}
for some  $h=(h^1,h^2,\cdots)$ and independent one-dimensional Wiener processes $W^k_t$ and
$
Y^k_t:=\int_{\bR^m}z \wt
{N}_k(t,dz).
$
Note that $Y^k_t$ are  independent $m$-dimensional pure jump L\'evy processes
with L\'e{}vy measure of $\nu^k$.

 Furthermore by considering $u-v$, where $v$ is the solution of
$$
dv=a(\omega,t)\Delta^{\alpha/2}v\,dt\,
+\sum_{i=k}^{\infty}h^k dW^k_t, \quad u(0)=0
$$
 from Theorem \ref{thm linear},  we find that without loss of generality we may also assume $h^k$'s are all zero.

By $<M,N>$ we denote the bracket of real-valued square integrable martingales $M$ and $N$. Also let $[M]$ denote the quadratic variation of $M$.
\begin{remark}
                  \label{remark 111}
(i) Note that, if $\wh{c}_{k,2}<\infty$, then  $Y^{k,i}=\int_{\bR^m}z^i\wt{N}
_k(t,dz)$ is a square integrable martingale for each
$k\geq 1$ and $i=1,\cdots,m$. Also  for any
  $\overline{\cP}$-measurable
 process $H=(H^1, \dots, H^m)\in L_2(\Omega\times [0,T], \bR^m)$ which has a predictable version
 $\bar{H}=(\bar{H}^1, \dots, \bar{H}^m)$,
 $$M^k_t:=\int_0^t H_s \cdot dY^k_s  =\sum_{i=1}^m\int_0^t \int_{\bR^m} H^{i}_s z^{i} \wt{N}
_k(ds,dz) = \sum_{i=1}^m\int_0^t \int_{\bR^m} \bar{H}^{i}_s z^{i} \wt{N}
_k(ds,dz)  $$
  is a square integrable
martingale with
$$
[M^k]_t=\sum_{i,j=1}^m \int^t_0\int_{\bR^m}H^iH^jz^iz^jN(ds,dz),
$$
$$
 E[M^k]_t=\sum_{i,j}(\int_{\bR^m}z^iz^j\nu^k(dz)) \E \int^t_0 H^i(s)H^j(s)ds\leq  \wh{c}^2m^2\E \int_0^t |H_s|^2 ds.
$$

 \noindent
(ii)
Suppose that  (\ref{main assump}) holds. Then, for any $1 \le j \le m$,
 $g^{\cdot, j}\in \bH^{\gamma}_p(T,\ell_2)$ and $\phi\in C^\infty_0
(\bR^d)$, the series of stochastic integral
$$\sum_{k=1}^{\infty}\sum_{j=1}^{m}\int^t_0(g^{k,j} (s, \cdot),\phi)dY^{k,j}_s
$$ defines a square
integrable martingale on $[0, T]$, which is right continuous with
left limits. Indeed, denote $M_n:=\sum_{k=1}^n\sum_{j=1}^m \int^t_0(g^{k,j} (s, \cdot),\phi)dY^{k,j}_s$, then the quadratic variation of $M_n$ is
 $$
[M_n]_t= \sum_{k=1}^n\sum_{i,j=1}^m  \int^t_0\int_{\bR^m}(g^{k,i},\phi)(g^{k,j}(s, \cdot),\phi) z^iz^kN^k(ds,dz),
$$
and
\begin{align*}
 &\E[M_n]_t=\sum_{k=1}^n\sum_{i,j=1}^m \E \int^t_0(g^{k,i},\phi)(g^{k,j}(s, \cdot),\phi)  \int_{\bR^m} z^i  z^j   \nu_k(dz)  ds
 \le c(m,\wh{c})\sum_{k=1}^n \sum_{j=1}^m\E\int^t_0(g^{k,j}(s, \cdot),\phi)^2 ds.
 \end{align*}
 Also, with $q:=p/(p-2)$,  for every
$1 \le j \le m$,
 \begin{align*}
& \sum_{k=1}^{\infty}\, \E \left[ \int^T_0(g^{k,j}(s, \cdot),\phi)^2ds \right]
= \sum_{k=1}^{\infty} \E \left[ \int_0^T
((1-\Delta)^{\gamma/2}g^{k,j}(s, \cdot),
(1-\Delta)^{-\gamma/2}\phi)^2 \, ds \right] \\
  \leq&
\|(1-\Delta)^{-\gamma/2}\phi\|_1 \, \E \left[ \int^T_0 \Big (
\sum_{k=1}^{\infty} |(1-\Delta)^{\gamma/2}g^{k,j}(s, \cdot)|^2, \,
|(1-\Delta)^{-\gamma/2}\phi| \Big) \, ds \right] \\
\leq& \|(1-\Delta)^{-\gamma/2}\phi\|_1 \,
\|(1-\Delta)^{-\gamma/2}\phi\|_q \, \E \left[ \int^T_0 \Big\|\,
 \sum_{k=1}^{\infty} |(1-\Delta)^{\gamma/2}g^{k,j}(s, \cdot)|^2    \Big\|_{p/2}\,
 ds \right]\\
\leq&  \|(1-\Delta)^{-\gamma/2}\phi\|_1 \,
\|(1-\Delta)^{-\gamma/2}\phi\|_q \, T^{1-\frac2{p}} \,
\|g^{\cdot, j}\|_{\bH^\gamma_p(T, l^2)}^2 <\infty.
\end{align*}
It follows that $[M_n]_t$ converges in probability uniformly on $[0,T]$ and this certainly proves the claim.
\end{remark}

\begin{defn}\label{D:2.5}
Write $u \in \cH^{\gamma+\alpha}_p(T)$ if $u\in
\bH^{\gamma+\alpha}_p(T), u(0)\in U^{\gamma+\alpha-\alpha/p}_p$, and
for some $f\in \bH^{\gamma}_p(T)$
$ h \in  \bH^{\gamma+\alpha/2}_p(T,\ell_2)$
and
 $ g^{\cdot,j} \in  \bH^{\gamma+\alpha/2}_p(T,\ell_2), 1 \le j \le m$
$$
du=f\,dt
+
\sum_{k=1}^{\infty}h^k dW^k_t +
\sum_{k=1}^{\infty}\sum_{j=1}^{m}g^{k,j}  d Y^{k,j}_t, \quad u(0)=u_0, \quad \hbox{for }  t\in [0,  T]
$$
 in the sense of distributions, that is, for any $\phi\in C^\infty_0 (\bR^d)$,
\begin{equation}
                \label{eqn 11.16}
(u(t, \cdot ),\phi)=(u(0, \cdot),\phi)+\int^t_0(f(s , \cdot),\phi)ds +\sum_{k=1}^{\infty} \int^t_0(h^k(s , \cdot),\phi)ds +
\sum_{k=1}^{\infty}\sum_{j=1}^m\int^t_0(g^{k,j}(s, \cdot),\phi)dY^{k,j}_s
\end{equation}
holds for all $t\leq T$ $a.s.$.
In this case we write
$$
\bD u :=f, \quad   \bS_c u :=(h^{1},\dots h^{k}, \dots), \quad
 \bS^{k,j}_d u :=g^{k,j}, \quad \bS^{\cdot,j}_d u :=(g^{1,j},\dots g^{k,j}, \dots)
$$
and define
$$
\|u\|_{\cH^{\gamma+\alpha}_p(T)}:=\|u\|_{\bH^{\gamma+\alpha}_p(T)}
 +\|\bD u\|_{\bH^{\gamma}_p(T)} +\|\bS_c u\|_{\bH^{\gamma+\alpha/2}_p(T,\ell_2)}+
\sum_{j=1}^m \|\bS^{\cdot,j}_d u\|_{\bH^{\gamma+\alpha/2}_p(T,\ell_2)}
+\|u(0)\|_{U^{\gamma+\alpha-\alpha/p}_p}.
$$
\end{defn}

To prove that $\cH^{\gamma+\alpha}_p(T)$ is a Banach space we need the following result,
which is an  infinite dimensional extension of Kunita's inequality (for example, see \cite[Theorem 4.4.23]{A}). In fact, if $m=1$ then the proof is given in \cite{C-K}.

\begin{lemma}
                      \label{littlewood}
Suppose
$1\le j \le m$,
$g^{\cdot, j}(\omega,t)=(g^{1, j},g^{2,j},\cdots)$'s are   $\ell_2$-valued predictable processes such that each $g^{k}=(g^{k,1}, \dots g^{k, m})$ is bounded.  Then, under the assumption \eqref{main assump},
\begin{eqnarray}
&&  \E\left[\left(\sum_{k=1}^{\infty}\int^t_0\int_{\bR^m} |g^k(s)|^2\,|z|^2
N_k
(s, dz)ds\right)^{p/2}\right] \nonumber \\
 &\leq&
c(p) \, \E \left[\left(\int^t_0\sum_{k=1}^{\infty}|g^k(s)|^2ds\right)^{p/2}+
\int^t_0\sum_{k=1}^{\infty}|g^k(s)|^p\, ds \right] \label{eqn 06.16.3}.
\end{eqnarray}
\end{lemma}

\pf   Due to  monotone convergence theorem we may assume $g^{k, j}=0$ for all $i>M$ and $1\le j \le m$. By monotone convergence theorem,
\begin{eqnarray*}
A&:=&\E \left[\left(\sum_{k=1}^{M}\int^t_0\int_{\bR^m}|g^k(s)|^2|z|^2
N_k(s,dz)ds\right)^{p/2}\right]\\
&=&\lim_{N\to \infty}\E
\left[ \left(\sum_{k=1}^{M}\int^t_0\int_{|z|\leq N}|g^k(s)|^2|z|^2
N_k(s,dz)ds\right)^{p/2}\right].
\end{eqnarray*}
Since
$(a+b)^{p/2}\leq c(p)(|a|^{p/2}+|b|^{p/2})$ and $\wt{N}
_k(s,dz):=
N_k(s,dz)-s\nu_k(dz)$,
\begin{eqnarray*}
A\leq  c(p)\lim_{N\to \infty} \E \left[ (J_{2,t})^{p/2} \right]
 +c(p) \E\left[ \left( \int^t_0\int_{\bR^m}
\sum_{k=1}^M |g^k(s)|^2\,|z|^2\nu_k(dz)ds\right)^{p/2}
\right]
 \end{eqnarray*}
where $J_{n,t}:=\sum_{k=1}^{M}\int^t_0\int_{|z|\leq N}|g^k(s)|^n|z|^n
\wt{N}_k(s,dz)ds$, which is a square integrable martingale becuase $g^k$ are bounded predictable  processes. 
By Burkholder-Davis-Gundy inequality (For example, see \cite[Theorem 48]{P}.)
\begin{align}
 &\E \left[ (J_{2,t})^{p/2} \right]\leq c(p) \E \left[[J_2 ]_t^{p/4} \right]=c(p)
 \E \left[\left(
\sum_{k=1}^{M}\int^t_0\int_{|z|\leq N}
 |g^k(s)|^4|z|^4 
N_k(s,dz)ds\right)^{p/4}\right] \label{eqn 6.1.1}\\
 &\le c(p)
\E\left[\left( \sum_{k=1}^M \,\,\sum_{0\leq s\leq t}|g^k (s)|^4 \,
|\Delta Y^k_s|^4\right)^{p/4}\right].\nonumber
 \end{align}
Recall that for any $q>1$, $(\sum |a_n|^q)^{1/q}\leq \sum |a_n|$.
Thus if $2< p \leq 4$, then
\begin{eqnarray*}
 \E \left[ (J_{2,t})^{p/2} \right]  \leq  c(p)\E \left[ \sum_{k=1}^M\,\,\sum_{0\leq
s\leq t}|g^k (s)|^p \, |\Delta Y^k_s|^p \right] \leq c(p,\wh{c})\E \left[ \int^t_0\sum_{k=1}^{M}|g^k (s)|^p
 \, ds \right] .
\end{eqnarray*}
If $4<p\leq 8$ then,  by  the relation  $\wt{N}_k(s,dz)=
N_k(s,dz)-s\nu_k(dz)$ and  Burkholder-Davis-Gundy inequality,
\begin{eqnarray*}
&& \E \left[ (J_{2,t})^{p/2} \right]  \leq c(p) \E \left[ \left(\sum_{k=1}^{M}\int^t_0\int_{|z|\leq N}|g^k(s)|^4 \, |z|^4 
N_k(s,dz)ds\right)^{p/4}\right]\\
&\leq& c(p)
\E \left[ (J_{8,t})^{p/8} \right]  +
c(p)\E\left[\left( \sum_{k=1}^{M}\int^t_0\int_{|z|\leq N}|g(s)|^4 \, |z|^4\nu_k(dz)ds\right)^{p/4}\right] \\
&\leq& c(p,\wh{c})\E \left[\left(\sum_{k=1}^{M}\int^t_0\int_{|z|\leq
N}|g^k (s)|^8|z|^8 
N_k(s,dz)ds\right)^{p/8} +\left(
\int^t_0\sum_{k=1}^{M}|g^k (s)|^4
ds \right)^{p/4}\right]\\
&\leq& c(p,\wh{c})\E \left[\left(\sum_{k=1}^{M}\int^t_0\int_{|z|\leq
N}|g^k (s)|^p|z|^p
N_k(s,dz)ds\right) +\left(
\int^t_0\sum_{k=1}^{M}|g^k (s)|^4
ds \right)^{p/4}\right]
\\
&\leq& c(p,\wh{c})\E\left[ \int^t_0\sum_{k=1}^{M}|g^k (s)|^p \, ds+ \left(\int^t_0\sum_{k=1}^{M}|g^k (s)|^4 ds
\right)^{p/4}\right].
\end{eqnarray*}
 Similarly, in general, for $p\in (2^{n-1}, 2^n]$,
 \begin{eqnarray*}
 A \leq c(p,\wh{c})  \, \sum_{j=1}^n \E \left[  \left(\int^t_0\sum_{k=1}^M|g^k (s)|^{2^j} ds\right)^{p2^{-j}} \right]  +c(p,\wh{c})\, \E \left[\int^t_0\sum_{k=1}^{M}|g^k(s, x)|^pds\right].
 \end{eqnarray*}
 Also since for each $2\leq q\leq
 p$,
$$
\left(\int^t_0\sum_{k= 1}^M |g^k (s)|^q
ds\right)^{1/q}\leq
\left(\left(\int^t_0\sum_{k=1}^M
 |g^k(s)|^2 ds\right)^{1/2}+\left(\int^t_0
 \sum_{k=1}^M |g^k(s)|^p ds\right)^{1/p}\right),
$$
we get
\begin{eqnarray}
                             \label{e: 12.6.4}
A \leq
c(p,\wh{c}) \, \E \left[\left(\int^t_0\sum_{k=1}^{\infty}|g^k(s)|^2ds\right)^{p/2}+ \int^t_0\sum_{k=1}^{\infty}|g^k(s)|^p\, ds \right].
\end{eqnarray}
Thus the lemma is proved.
\qed

\begin{thm}
                        \label{theorem banach}
            Suppose that  \eqref{main assump} holds.
For any $p\in [2, \infty)$ and $\gamma\in \bR$,
$\cH^{\gamma+\alpha}_p(T)$ is a Banach space
with norm $\| \cdot\|_{\cH^{\gamma+\alpha}_p(T)}$. Moreover, there
is a constant $c=c(d,p,T)>0$ such that for every $u\in
\cH^{\gamma+2}_p(T)$ and $0 < t \le T$,
\begin{equation}
                                         \label{e:2.6}
\E \left[ \sup_{s\leq t}\|u(s, \cdot )\|^p_{H^{\gamma}_p} \right] \leq c(p,d, T)
\left( \|\bD u\|^p_{\bH^{\gamma}_p(t)}+\|\bS_c u\|_{\bH^{\gamma}_p(t,\ell_2)}
+\sum_{j=1}^m \|\bS^{\cdot,j}_d   u\|^p_{\bH^{\gamma }_p(t,\ell_2)}
+\|u_0\|^p_{U^{\gamma}_p} \right) .
\end{equation}
\end{thm}

\pf
 By Theorem \ref{thm banach conti} and the reasons explained just before Remark \ref{remark 111}, without loss of generality we assume that
 $Y^k_t=\int_{\bR^m}z\wt{N}
_k(t,dz)$. Moreover, due to  Remark \ref{remark 11.22} it suffices to prove the
theorem only for $\gamma=0$.  First we prove (\ref{e:2.6}). Let
$du=fdt +\sum_{k=1}^{\infty} g^k \cdot dY^k_t$ with $u(0)=u_0$.

For a moment, we assume that $g^{k,j}=0$ for all
$k\geq N_0, 1\le j \le m$ and  $g^{k,j}$ is of the type
\begin{equation}
                               \label{e:12.5}
g^{k,j}(t,x)=\sum_{i=0}^{m_k} I_{(\tau^{k,j}_{i},\tau^{k,j}_{i+1}]}(t)g^{k_i, j}(x),
\end{equation}
where $\tau^{k,j}_i$ are bounded stopping times and $g^{k_i, j}\in
C^{\infty}_0(\bR^d)$. Define
$$
v(t,x)=\sum_{k=1}^{N_0} \int^t_0 g^k (s, x) \cdot dY^k_s.
$$
Then by Burkholder-Davis-Gundy inequality and Lemma \ref{littlewood},
\begin{eqnarray*}
&&\E \left[ \sup_{s\leq t}|v(s,x)|^p \right] =\E \left[ \sup_{s\leq t}\left|\sum_{k=1}^{N_0} \sum_{j=1}^{m}\int^t_0 \int_{\bR^m} g^{k,j} (s, x) z^j
\wt{N}
_k(s,dz)ds
        \right|^p \right]\\
        &\leq&  c(p) \E
\left[\left(\sum_{k=1}^{\infty}\sum_{i, j=1}^{m}\int^t_0\int_{\bR^m}  g^{k,i}(s,x)g^{k,j}(s,x)z^iz^j
N_k(s,dz)ds\right)^{p/2}\right]\\
&\leq&  c(p,\wh{c}) \E
\left[\left(\sum_{k=1}^{\infty}\int^t_0\int_{\bR^m} |g^k(s,x)|^2|z|^2
N_k(s,dz)ds\right)^{p/2}\right]\\
&\leq& c(p,\wh{c})\, \E\left[ \left(\int^t_0\sum_{k=1}^\infty |g^k (s, x)|^2ds \right)^{p/2}+ \int^t_0\sum_{k=1}^{\infty}|g^k(s, x)|^pds\right] .
 \end{eqnarray*}
Since $\sum_n |a_n|^p\leq (\sum_n |a_n|^2)^{p/2}$ and $(\int^t_0 |f|ds)^p\leq t^{p-1}\int^t_0|f|^pds$,  by integrating over $\bR^d$ we get that
for every $t \le T$
\begin{equation}
                     \label{eqn 11.22.1}
\E \left[\sup_{s\leq t}\|v\|^p_p \right]\leq c(T,p) \sum_{j=1}^m
\|g^{\cdot, j}\|^p_{\bL_p(t,\ell_2)}:=c(T,p)\, \sum_{j=1}^m\E \int^t_0 \int_{\bR^d} |g^{\cdot, j}|^p_{\ell_2}\,dxds.
\end{equation}
Next we prove \eqref{eqn 11.22.1}  for general $g^{\cdot, j} \in \bL_p(T,\ell_2)$. By Theorem 3.10 in \cite{Kr99}, we can
 take a sequence $g^{\cdot, j}_n\in \bL_p(T,\ell_2)$ so that for each fixed $n$,
$g^{k,j}_n=0$ for all large $k$ and each $g^{k,j}_n$ is of of the type
(\ref{e:12.5}), and $g^{\cdot, j}_n \to g^{\cdot, j}$ in $\bL_p(T, \ell_2)$ as $n\to
\infty$. Define $v_n(t,x)=\sum_{k=1}^{\infty}\sum_{j=1}^{m}\int^t_0 g^{k,j}_n dY^{k,j}_t$, then for every $t \le T$
$$
\E \left[ \sup_{s\leq t}\|v_n\|^p_p\right] \leq c(T,p)
\sum_{j=1}^{m}\|g^{\cdot, j}_n\|^p_{\bL_p(t,\ell_2)}, \quad \E \left[\sup_{s\leq
t}\|v_{n_1}-v_{n_2}\|^p_p\right] \leq c(T,p)
\sum_{j=1}^{m} \|g^{\cdot, j}_{n_1}-g^{\cdot, j}_{n_2}\|^p_{\bL_p(t,\ell_2)}.
$$
Thus
\eqref{eqn 11.22.1} follows by taking $n\to \infty$.
Now note that
$$
d(u-v)=f dt  \quad \hbox{ with} \quad (u-v)(0)=u_0.
$$
Thus it is easy to check that
$$
\E \left[ \sup_{s\leq t}\|u-v\|^p_p \right] \leq N\E \left[
\|u_0\|^p_p\right] +N\E \left[ \int^t_0\|f(s, \cdot) \|^p_p \,
ds\right].
$$
Consequently,
$$
\E \left[ \sup_{s\leq t}\|u\|^p_p \right] \leq
 N
\|f\|^p_{\bH^0_p(t)}+c\sum_{j=1}^{m} \|g^{\cdot,j}\|^p_{\bL_p(t,\ell_2)}+N\E\|u_0\|^p_{p}.
$$
The completeness of the space $\cH^{\alpha}_p(T)$ easily follows from
\eqref{e:2.6}.  Indeed, let $\{u_n: n=1,2,\cdots\}$  be a Cauchy sequence in $\cH^{\alpha}_p(T)$. Then $\{u_n\}$, $\{\bD u_n\}$, $\{\bS^{\cdot,j}_d u_n\}$ and $\{u_n(0)\}$ are Cauchy sequences in $\bH^{\alpha}_p, \bL_{p}(T), \bH^{\alpha/2}_p(T,\ell_2)$ and $U^{\alpha/2-\alpha/p}_p$ respectively. Thus there exist
$u\in \bH^{\alpha}_p(T)$, $f\in \bL_{p}(T), g^{\cdot,j}\in \bH^{\alpha/2}_p(T,\ell_2)$ and $u_0\in U^{\alpha/2-\alpha/p}_p$ so that
$u_n, \bD u_n, \bS^{\cdot,j}_d u_n, u_n(0)$ converge to $u, f,g^{\cdot,j},u_0$ respectively, that is,
$$
\|u_n-u\|_{\bH^{\alpha}_p(T)}+\|\bD u_n -f\|_{\bL_p(T)}
+\|\bS^{\cdot,j} u_n -g^j\|_{\bH^{\alpha/2}_p(T,\ell_2)}+\|u_n(0)-u_0\|_{U^{\alpha/2-\alpha/p}_p}\to 0
$$
as $n \to \infty$. Thus to prove $u\in \cH^{\alpha}_p(T)$ and $u_n \to u$ in $\cH^{\alpha}_p(T)$, we only need to show that for any $\phi\in C^{\infty}_0(\bR^d)$, the equality
\begin{equation}
                           \label{equality}
(u(t, \cdot),\phi)=(u_0,\phi)+\int^t_0(f(s, \cdot),\phi)ds+\sum_{k=1}^{\infty} \sum_{j=1}^m\int^t_0(g^{k,j}(s, \cdot),\phi)dY^{k,j}_s
\end{equation}
holds for all $t\leq T$ (a.s.).  Taking the limit from
$$
(u_n(t, \cdot),\phi)=(u_n(0),\phi)+\int^t_0(\bD u_n(s, \cdot),\phi)ds+ \sum_{k=1}^{\infty}\sum_{j=1}^m\int^t_0(\bS^{k,j}_d u_n(s, \cdot),\phi)dY^{k,j}_s
$$
and using the argument used in Remark \ref{remark 111}(ii) one can show that  (\ref{equality}) holds in $\Omega \times [0,T]$ (a.e.).  Also using the inequality (see \eqref{e:2.6})
$$
\E \left[ \sup_{t\leq T}\|u_n(\cdot, t)-u_m(\cdot, t)\|^p_{L_p} \right]\leq N \|u_n-u_m\|_{\cH^{\alpha}_p(T)}
$$
and taking $m \to \infty$, one finds that $(u(t, \cdot),\phi)$ is right continuous with left limits, and consequently (\ref{equality}) holds for all $t\leq T$ (a.s.). The theorem is proved.
\qed

\begin{lemma}
                   \label{lemma 3.2}
Let $p\in (2,\infty)$, $t>0$ and $f\in L_p([0,t \, ]\times \bR^d)$. Then for
any $\varepsilon>\alpha(1/2-1/p)$,
\begin{equation}
                       \label{e:12.1.1}
\int_{\bR^d}\int^{t}_0\int^s_0|\partial^{\alpha/2}_xT_{s-r}f(r,x)|^p\,drds\,dx\leq
c\int^{t}_0\|f(s, \cdot)\|^p_{H^{\varepsilon}_p}\, ds,
\end{equation}
where $c=c(d,p,\alpha,\varepsilon)$ is independent of $t$.
\end{lemma}

\pf  Note that we may  assume  $\alpha(1/2-1/p)<\varepsilon<\alpha/2$. Let $q>p$ be
chosen so that
$$
\frac{1}{p}=(1-\frac{2\varepsilon}{\alpha})\times \frac{1}{2}+\frac{2\varepsilon}{\alpha}\times \frac{1}{q}.
$$
Such choice of $q$ is possible since $1/p>(1-\frac{2\varepsilon}{\alpha})\times \frac{1}{2}$.
We will use  an interpolation theorem.  First, note that
 $$
 \varepsilon=(1-\frac{2\varepsilon}{\alpha})\times 0+\frac{2\varepsilon}{\alpha}\times \frac{\alpha}{2}.
 $$
 Define
 an operator $\cA$ by
$$
\cA f(s,r,x)=\begin{cases} \partial^{\alpha/2} T_{s-r}f  \quad &\hbox{if } r<s, \\
 0 &\hbox{otherwise}. \end{cases}
 $$
 Then, due to (\ref{eqn ildoo}) and the inequality
 $\|T_{s-r}\partial^{\alpha/2}f\|_q\leq \|\partial^{\alpha/2}f\|_q\leq \|f\|_{H^{\alpha/2}_q}$,
the linear mappings
 $$
\cA: L_2([0,t],L_2(\bR^d)) \to L_2([0,t]\times [0,t]\times \bR^d)
$$
 and
$$
\cA: L_q([0,t],H^{\alpha/2}_q) \to L_q([0,t]\times [0,t]\times \bR^d)
$$
are bounded and their norms are independent of $t$.
 It follows from the interpolation theory (see, for
instance, \cite[Theorem 5.1.2]{BL}) that the operator
$$
\cA: L_p([0, t],H^{\varepsilon}_p(\bR^d)) \to L_p([0,t]\times
[0, t]\times \bR^d)
$$
is bounded and its norm is independent of $t$.
 The lemma is proved.
 \qed

\begin{thm}
                   \label{thm main levy}
Fix a constant $\varepsilon_1$ so that $\varepsilon_1=0$ if $p=2$, and $\varepsilon_1>\alpha(1/2-1/p)$ if $p>2$.
suppose \eqref{main assump} holds.
Then for any $f\in \cH^{\gamma}_p(T), h\in \bH^{\gamma+\alpha/2}_{p}(T,\ell_2), g^{\cdot, j}\in \cH^{\gamma+\alpha/2+\varepsilon_1}_p(T,\ell_2), 1 \le j \le m$ and $u_0\in U^{\gamma+\alpha/2-\alpha/p}_p$,  equation (\ref{eqn main levy}) has a unique solution $u$ in $\cH^{\gamma+\alpha}_p(T)$, and for this solution
\begin{equation}
                     \label{eqn linear levy}
\|u\|_{\cH^{\gamma+\alpha}_p(t)}\leq c(p,T,\delta)\left(\|f\|_{\bH^{\gamma}_p(t)}+\|h\|_{\bH^{\gamma+\alpha/2}_p(t,\ell_2)}
+\sum_{j=1}^m\|g^{\cdot, j}\|_{\bH^{\gamma+\alpha/2+\varepsilon_1}_p(t,\ell_2)}+\|u_0\|_{U^{\gamma+\alpha-\alpha/p}_p}\right)
\end{equation}
 for every $t \le T$.
\end{thm}

 \pf
 As explained before, without loss of generality we assume that
$h^i$'s are all zeros.

{\bf{Step 1}}. As in the proof of Theorem \ref{thm linear},  we only need to prove the theorem
 for a particular $\gamma=\gamma_0$.

{\bf Step 2}. We assume $a(\omega)=1$ and  prove the theorem for the equation:
\begin{equation}
                        \label{e: heat eqn}
  du=\Delta^{\alpha/2} u dt + \sum_{k=1}^{\infty} g^k dY^k_t, \quad u(0)=0.
  \end{equation}
By the result of Step 1, we may assume that $\gamma=-\alpha/2$. The uniqueness is obvious and we only prove the existence and the estimate \eqref{eqn linear levy}.
Considering approximation arguments, for a moment, we assume that
$g^{k,j}=0$ for all $k > N_0$ and $1 \le j \le m$ and that
$$
g^{k, j} (t,x)=\sum_{i=0}^{m_k}I_{(\tau^{k,j}_i,\tau^{k,j}_{i+1}]}(t)g^{k_i,j}(x),
$$
where $\tau^{k,j}_i$ are bounded stopping times and $g^{k_i,j}(x)\in
 C^{\infty}_0 (\bR^d)$. Define
$$
v(t,x):=\sum_{k=1}^{N_0} \int^t_0
 g^k(s, x) \cdot dY^k_s
=\sum_{k=1}^{N_0}\sum_{i=1}^{m_k} \sum_{j=1}^{m}g^{k_i, j}(x)(Y^{k,j}_{t\wedge
\tau^k_{i+1}}-Y^{k,j}_{t\wedge \tau^k_{i}})
$$
and
\begin{equation}
                    \label{eqn 03-22-11}
u(t,x):=v(t,x)+\int^t_0 \Delta^{\alpha/2}T_{t-s} v \,
ds=v(t,x)+\int^t_0 T_{t-s}\Delta^{\alpha/2} v \, ds.
\end{equation}
 Now we remember from the proof of Theorem \ref{thm linear} that if  functions $h_1=h_1(t,x)$ and $h_2=h_2(x)$ are  sufficiently smooth, then
$$
w_1(t,x):=\int^t_0T_{t-s}h_1(s) ds, \quad w_2(t,x)=T_t h_2
$$
solve
$$
dw_1=(\Delta^{\alpha/2}w + h_1)\,dt, \quad w_1(0)=0,
$$
$$
dw_2=\Delta^{\alpha/2}w_2\,dt, \quad w_2(0)=h_2.
$$
Therefore we have $d(u-v)=(\Delta^{\alpha/2}(u-v)+\Delta^{\alpha/2}
v)dt=\Delta^{\alpha/2} u dt$, and
$$
du=\Delta^{\alpha/2} u dt +dv=\Delta^{\alpha/2} u dt +\sum_{k=1}^{N_0} g^k  \cdot dY^k_t.
$$
Let $T_{t-s}g^k (r, x)=(T_{t-s}g^{k,1} (r, x), \dots T_{t-s}g^{k, m} (r, x))$.
By (\ref{eqn 03-22-11}) and stochastic Fubini theorem (\cite[Theorem 64]{P}), almost
surely,
\begin{eqnarray} \label{e:fger2}
u(t,x)&=&v(t,x)+\sum_{k=1}^{N_0}\int^t_0 \int^s_0
\Delta^{\alpha/2}T_{t-s}
g^k (r, x) \cdot dY^k_r ds \nonumber \\
&=&v(t,x)-\sum_{k=1}^{N_0} \sum_{j=1}^m\int^t_0\int^t_r\frac{\partial}{\partial
s}T_{t-s}g^{k,j} (r, x) ds dY^{k,j}_r \nonumber\\
 &=&\sum_{k=1}^{N_0} \int^t_0 T_{t-s}g^k (s,x) \cdot dY^k_s.
\end{eqnarray}
Hence,
$$
\partial^{\alpha/2}_xu(t,x)=\sum_{k=1}^{N_0}\int^t_0
\partial^{\alpha/2}_xT_{t-s}g^k ( s,x)  \cdot dY^k_s=
\sum_{k=1}^{N_0}\sum_{j=1}^{m}\int^t_0
\partial^{\alpha/2}_xT_{t-s}g^{k,j}( s,x)  dY^{k,j}_s
.
$$
By
 Burkholder-Davis-Gundy's
inequality and Lemma \ref{littlewood}, we have for every $0<t \le T$
\begin{align*}
&\E \left[  |\partial^{\alpha/2}_xu(t,x)|^p \right]
\leq c(p)\E\left[
\left(\sum_{k=1}^{N_0}\sum_{i,j=1}^{m} \int^t_0\int_{\bR^m}
\partial^{\alpha/2}_xT_{t-s}g^{k,i}( s,x) \partial^{\alpha/2}_xT_{t-s}g^{k,j}( s,x)  z^i z^j    N^k(dz,ds)\right)^{p/2}\right]  \\&\leq\, c(p)\E\left[
\left(\sum_{k=1}^{N_0}\int^t_0\int_{\bR^m}
|\partial^{\alpha/2}_xT_{t-s}g^k( s,x)|^2|z|^2N^k(dz,ds)\right)^{p/2}\right]  \\
&\leq  c(p)\, \E\left[ \left(
\int^t_0\sum_{k=1}^{\infty}|\partial^{\alpha/2}_xT_{t-s}g^k( s,x)|^2ds\right)^{p/2}
+
\int^t_0\sum_{k=1}^{\infty}|\partial^{\alpha/2}_xT_{t-s}g^k( s,x)|^pds\right].
 \end{align*}
 By (\ref{eqn ildoo}), Lemma \ref{lemma 3.2} and the inequality
$\sum_{k=1}^{\infty} |a_k|^p\leq (\sum_{k=1}^{\infty} |a_n|^2)^{p/2}$,
\begin{equation}
                               \label{e:12.2.31}
\E \left[ \int^t_0\|\partial^{\alpha/2}_xu (s, \cdot) \|^p_{p} \, ds \right]\leq c(p, \alpha)  \sum_{j=1}^m\E
\left[ \int^t_0\|g^{\cdot, j}(s, \cdot)\|^p_{H^{\varepsilon_1}_p(\ell_2)} \, dt\right].
\end{equation}
Similarly from \eqref{e:fger2} we also get,  for every $0<t \le T$,
\begin{equation}
                                 \label{eqn 8.6.51}
\E \left[ |u(t,x)|^p \right] \leq c(p)\, \E\left[ \left(
\int^t_0\sum_{k=1}^{N_0}|T_{t-s}g^k (s,x) |^2ds\right)^{p/2} +
\int^t_0\sum_{k=1}^{N_0}|T_{t-s}g^k(s,x)|^p\, ds\right].
 \end{equation}
By the same argument which leads to \eqref{e:hjkrty}, we see that the right side of \eqref{eqn 8.6.51} is finite.
 Thus we proved $\partial^{\alpha/2}_x u, u\in \bL_{p}(T)$, and  hence $u\in \cH^{\alpha/2}_p(T)$. As in (\ref{e1}) and (\ref{e2}),
\begin{eqnarray*}
 \|u\|^p_{\cH^{\alpha/2}_p(t)}&\leq& c(p)\left(\|u\|^p_{\bH^{\alpha/2}_p(t)}+\|\Delta^{\alpha/2}u\|^p_{\bH^{-\alpha/2}_p(t)}+\|g\|^p_{\bL_p(t,\ell_2)}\right)\\
 &\leq& c\left(\|u\|^p_{\bH^{-\alpha/2}_p(t)}+\|\partial^{\alpha/2}_xu\|^p_{\bL_p(t)}+\|g\|^p_{\bL_p(t,\ell_2)}\right)\\
 &\leq& c(p,T,\alpha)\left(\|u\|^p_{\bH^{-\alpha/2}_p(t)} +
\|g\|^p_{\bH^{\varepsilon_1}_p(t,\ell_2)} \right)\\
&\leq& c(p,T,\alpha) \int_0^t\|u\|^p_{\cH^{\alpha/2}_p(s)} ds +  c(p,T,\alpha)
\|g\|^p_{\bH^{\varepsilon_1}_p(T,\ell_2)}.
 \end{eqnarray*}
 Finally,  Gronwall  leads to \eqref{eqn linear levy}.
Once one has a unique solvability of equation (\ref{e: heat eqn}) and estimate
(\ref{eqn linear levy}) for sufficiently smooth $g$, we repeat the same approximation argument  used   in  the
 Step 2 of  the
proof of Theorem \ref{thm linear}.

 {\bf Step 3}. Now,   we follow Step 3--Step 5 of  the
proof of Theorem \ref{thm linear}   word for word except obvious changes from $W^k_t$ to $Y^k_t$. The theorem is proved.
\qed

Finally we consider the nonlinear equation
\begin{equation}
                                \label{eqn nonlinear1}
du=\left(a(\omega,t)\Delta^{\alpha/2}u+f(u)\right)\,dt
+
\sum_{k=1}^{\infty}h^k(u) dW^i_t +
\sum_{k=1}^{\infty}\sum_{j=1}^{m}g^{k,j} (u) d Y^{k,j}_t, \quad u(0)=u_0,
\end{equation}
where $f(u)=f(\omega,t,x,u)$,
$h^k(u)=h^k(\omega,t,x,u)$,
$g^k(u)=(g^{k,1}(\omega,t,x,u), \dots, g^{k,m}(\omega,t,x,u))$, $W_t$ are independent $1$-dimensional Wiener processes and
$
Y^k_t:=\int_{\bR^m}z \wt
{N}_k(t,dz)
$
are  independent $m$-dimensional pure jump L\'evy processes with L\'evy measure $\nu_k$.

\begin{assumption}
                 \label{ass nonlinear1}
                 Fix a constant $\varepsilon_1$ so that $\varepsilon_1=0$ if $p=2$, and $\varepsilon_1>\alpha(1/2-1/p)$ if $p>2$.
Assume  that        $f(0) \in \bH^{\gamma}_p(T)$, $g^{\cdot, j}(0) \in \bH^{\gamma+\alpha/2+\varepsilon_1}_p(T, \ell_2)$ and $h(0) \in \bH^{\gamma+\alpha/2}_p(T, \ell_2).$
 Moreover, for any $\varepsilon>0$, there exists a constant $K_{\varepsilon}$ so that
  for any $u=u(x),v=v(x) \in H^{\gamma+\alpha}_p$ and
  $t, \omega$
  we have
  \begin{align}
                       \label{eqn 03.18.911}
 & \|f(t, \cdot, u(\cdot))-f(t, \cdot, v(\cdot))\|_{H^{\gamma}_p} + \sum_{j=1}^m\|g^{\cdot, j}(t, \cdot, u(\cdot))-g^{\cdot, j}(t, \cdot, v(\cdot))\|_{H^{\gamma+\alpha/2+\varepsilon_1}_p(\ell_2)}\nonumber \\ & +\|h(t, \cdot, u(\cdot))-h(t, \cdot, v(\cdot))\|_{H^{\gamma+\alpha/2}_p(\ell_2)}
 \leq \varepsilon \|u-v\|_{H^{\gamma+\alpha}_p}+K(\varepsilon)\|u-v\|_{H^{\gamma}_p}.
  \end{align}
\end{assumption}

\begin{example} \rm
Recall that the space  $B^{r}$ is defined in \eqref{e:2.11}.
Fix $\kappa_0=\kappa_0(\gamma)\geq 0$ so that $\kappa_0>0$ if $\gamma$ is not integer.
Consider
\begin{align*}
f(u)&=b(\omega,t,x)\Delta^{\beta_1/2}u+ \sum_{i=1}^d c^i(\omega,t,x)u_{x^i}I_{\alpha>1}+d(\omega,t,x)u+f_0,\\
 h^k(u)&=\eta^{k}(\omega,t,x)\Delta^{
 \beta_2/2}u+l^k
(\omega,t,x)u+h^k_0,\\
 g^{k,j}(u)&=\sigma^{k,j}(\omega,t,x)\Delta^{
 \beta^j_3/2}u+v^{k,j}(\omega,t,x)u+g^{k,j}_0, \quad j=1, \dots m.
 \end{align*}
 Here $\beta_1<\alpha$, $\beta_2<\alpha/2$, $\beta^j_3<\alpha/2-\varepsilon_1$ and $f_0 \in \bH^{\gamma}_p(T)$, $h_0\in \bH^{\gamma+\alpha/2}_p(T,\ell_2)$, $g^{\cdot, j}_0\in \bH^{\gamma+\alpha/2+\eps_1}_p(T,\ell_2)$.
 Assume for each $\omega,t,i,j$,
 \begin{align*}
 &|b|_{B^{|\gamma|+\kappa_0}}+|c^i|_{B^{|\gamma|+\kappa_0}}+|d|_{B^{|\gamma|+\kappa_0}}
 +|\eta|_{B^{|\gamma|+\alpha/2+\kappa_0}}+|l|_{B^{|\gamma|+\alpha/2+\kappa_0}}\\
 &+
  |\sigma^{\cdot, j} |_{B^{|\gamma|+\alpha/2+\varepsilon_1+\kappa_0}}+|v^{\cdot, j}|_{B^{|\gamma|+\alpha/2+\varepsilon_1+\kappa_0}}\leq K < \infty.
 \end{align*}
Then the calculus in Example \ref{ex 1}  shows that (\ref{eqn 03.18.911}) holds.
\end{example}

Here is the main result of this section.

\begin{thm}\label{main thm}
Suppose \eqref{main assump} and Assumptions \ref{ass 1} and \ref{ass nonlinear1} hold.
Then
the  equation (\ref{eqn nonlinear1}) has a unique solution $u\in \cH^{\gamma+\alpha}_p(T)$, and for this solution we have
\begin{equation}
                          \label{eqn 6.02.2}
  \|u\|_{\cH^{\gamma+\alpha}_p(t)}\leq c\left(\|f(0)\|_{\bH^{\gamma}_p(t)}
  +\|h(0)\|_{\bH^{\gamma+\alpha/2}_p(t,\ell_2)}
  +  \sum_{j=1}^m \|g^{\cdot, j}(0)\|_{\bH^{\gamma+\alpha/2+\varepsilon_1}_p(t,\ell_2)}
  +\|u_0\|_{U^{\gamma+\alpha/2-\alpha/p}_p}\right),
\end{equation}
 for every $t \le T$, where $c=c(p,T,\delta)$.
 \end{thm}
\pf
As we mentioned in the previous section,
our proof is a repetition of    that of Theorem 6.4 in \cite{Kr99}.
By Theorem \ref{thm main levy}, for any $u\in \cH^{\gamma+\alpha}_p(T)$ with initial data $u_0$ we can define $v=\cR u$  as  the solution of
\begin{align*}
dv =\left(a(\omega,t)\Delta^{\alpha/2} v(t, x)+f(t, x, u)\right)dt
+
\sum_{i=1}^{\infty} h^i(t, x, u) dW^i_t+
\sum_{k=1}^{\infty} \sum_{j=1}^{m} g^{k,j} (t, x, u) d Y^{k,j}_t,  \quad  v(0)=u_0.
\end{align*}
Then for any $u,v$ initial data $u_0$,  we have $(\cR u-\cR v)(0, x)=0$ and
\begin{align*}
d(\cR u-\cR v)
=&\left(a(\omega,t)\Delta^{\alpha/2} (\cR u- \cR v)+(f(t, x, u)-f(t, x, v))\right)dt
\\
&+
\sum_{i=1}^{\infty}\int_0^t (h^i(t, x, u)-    h^i(t, x, u)) dW^i_t +
\sum_{k=1}^{\infty} \sum_{j=1}^{m}\int_0^t  ( g^{k,j} (t, x, u) -g^{k,j} (t, x, v) )d Y^{k,j}_t.
\end{align*}
By Theorems \ref{thm linear} and Assumption \ref{ass nonlinear1}, for every $t \in (0, T]$,
\begin{align*}
&\|\cR u -\cR v\|^p_{\cH^{\gamma+\alpha}_p(t)}\\
&\leq c(p,T,\delta)\Big(\|f(u)-f(v)\|^p_{\bH^{\gamma}_p(t)}
+\|h(u)-h(v) \|^p_{\bH^{\gamma+\alpha/2}_p(t,\ell_2)}
+\sum_{j=1}^m\|g^{\cdot, j}(u)-g^{\cdot, j}(v)\|^p_{\bH^{\gamma+\alpha/2+\varepsilon_1}_p(t,\ell_2)}\Big)\\
&\le
\varepsilon^p  c(p,T,\delta)\|u -v\|^p_{\cH^{\gamma+\alpha}_p(t)}+K(\varepsilon) c(p,T,\delta)\int_0^t \E \|u(s, \cdot) -v(s, \cdot)\|^p_{H^{\gamma}_p} ds\\
&\le
\theta\|u -v\|^p_{\cH^{\gamma+\alpha}_p(t)}+ N\int_0^t  \|u -v\|^p_{\cH^{\gamma+\alpha}_p(s)} ds
\end{align*}
 where $\theta:=\varepsilon^p  c(p,T,\delta)$ and $N=c(p,T,\delta,\varepsilon)$.  Denote $\cR^{n+1}u:=\cR (\cR^n u)$. Then by induction, for every $t \in (0, T]$
 \begin{eqnarray*}
 &&\|\cR^n u -\cR^n v\|^p_{\cH^{\gamma+\alpha}_p(t)} \\
 &\le &\theta^n \|u -v\|^p_{\cH^{\gamma+\alpha}_p(t)} +
 \sum_{k=1}^n {{n} \choose{k}} \theta^{n-k} N^k \int_0^t \frac{(t-s)^{k-1}}{(k-1) !} \|u -v\|^p_{\cH^{\gamma+\alpha}_p(s)} ds
  \end{eqnarray*}
Therefore,
\begin{eqnarray*}
\|\cR^n u -\cR^n v\|^p_{\cH^{\gamma+\alpha}_p(T)}
& \le &
 \theta^{n}\sum_{k=0}^n {{n} \choose{k}}   \frac{(NT/\theta)^{k}}{k !} \|u -v\|^p_{\cH^{\gamma+\alpha}_p(T)}\\
 & \le &
 (2\theta)^{n} \left(\sup_{ k \ge 0}  \frac{(NT/\theta)^{k}}{k !} \right) \|u -v\|^p_{\cH^{\gamma+\alpha}_p(T)}.
  \end{eqnarray*}
  Choose $\eps>0$  so that $2\theta<1/2$, and then fix $n$ large enough so  that $ (2\theta)^{n} \left(\sup_{ k \ge 0}  \frac{(NT/\theta)^{k}}{k !} \right)<1/2$. Then
  $\bar{\cR}:=\cR^n$  is a contraction in $\cH^{\gamma+\alpha}_p(T)$ and obviously the unique fixed point $u$ under this map becomes the unique solution of
 (\ref{eqn nonlinear1}).
   Moreover, the estimate  \eqref{eqn 6.02.2} also easily from  Assumption \ref{ass nonlinear1}, Theorems \ref{thm main levy} and \ref{theorem banach}. We leave the details to the readers as an exercise.
\qed

\section{Application and Extension}
                                      \label{section 4}

First, we consider  equations with the random fractional Laplacian driven by (L\'evy) space-time white noise;
Let $d=1$ and consider the equation
\begin{equation}
                       \label{space-time}
du=(a(\omega,t) \Delta^{\alpha/2} u(t,x)+f(\omega,t,x,u(t,x)))dt+ \xi(\omega,t,x)h(\omega,t,x,u(t,x))d\cZ_t
\end{equation}
where $\cZ_t$ is a cylindrical L\'evy process  on $L_2(\bR)$, that is $\cZ_t$ has an expansion of the form
$$
\cZ_t=\sum_{k=1}^{\infty} \eta^k(x) Z^k_t
$$
where  $\{\eta^k:k=1,2,\dots\}$ is an orthonormal basis in $L_2$ and  $Z^k_t$ are
i.i.d. one-dimensional $\cF_t$-adapted L\'evy processes (see \cite{HOUZ} for the details). Using this expansion we can rewrite (\ref{space-time}) as follows :
\begin{equation}
                       \label{space-time white}
du=(a(\omega,t) \Delta^{\alpha/2} u+f(u))dt+ \sum_{k=1}^{\infty}g^k(u)d Z^k_t,
\end{equation}
where $g^k(u):=\xi(\omega,t,x) h(\omega,t,x,u(t,x))\eta^k(x)$.

Let $\gamma,p,s,r$ be constants satisfying
\begin{equation}
                 \label{eqn 06.04.1}
0>\gamma +\alpha /2>-1, \quad p\geq 2r\geq 2, \quad 1 \le r < (2\gamma +\alpha +2)^{-1}, \quad  s^{-1}+r^{-1}=1 \,\,(1\le s \le \infty).
\end{equation}
Define
$$
R_{\gamma}(x):=|x|^{-(\gamma +\alpha/2+1)}\int^{\infty}_0 t^{-(\gamma+\alpha/2+3)/2} e^{-tx^2-1/(4t)}dt.
$$
It is known that there exists a constant $c>0$ so that $c R_{\gamma}(x)$ is the kernel of the operator $(1-\Delta)^{(\gamma+\alpha/2)/2}$, that is
 $(1-\Delta)^{(\gamma+\alpha/2)/2}f=(cR_{\gamma}*f)(x)$.

\begin{assumption}
                     \label{ass 06.06.1}
    (i) For each $x$, $\xi=\xi(\omega,t,x)$ is predictable, and  $\|\xi(\omega,t,\cdot)\|_{L_{2s}}\leq K$ for each $\omega,t$.

    (ii) For each $x,u$, the processes $f(\omega,t,x,u), h(\omega,t,x,u)$ are predictable, and
    $$
    |f(\omega,t,x,u)-f(\omega,t,x,v)|\leq K|u-v|, \quad   |h(\omega,t,x,u)-h(\omega,t,x,v)|\leq K|u-v|.
    $$
    \end{assumption}

By following the  arguments  in the proof of \cite[Lemma 8.4]{Kr99}, we get the following

\begin{lemma}\label{l last1}
Let (\ref{eqn 06.04.1}) hold. Take some functions $h_0=h_0(x) \in L_p(\bR)$, $\xi_0=\xi_0(x) \in L_{2s}(\bR)$, and set $g^k_0= \xi_0 h_0 \eta^k$. Then $g_0= \{g^k_0\} \in H^{\gamma+\alpha/2}_p (\ell_2)$ and
$$
\|g_0\|_{H^{\gamma+\alpha/2}_p(\ell_2)} = \|\overline{h}_{0,\gamma}\|_p \le  N \|\xi_0\|_{2s} \|h_0\|_p,
$$
where $N=\|R_{\gamma}\|_{2r}<\infty$ and
$$
\bar{h}_{0,\gamma}(x):=\left(\int_{\bR}R_{\gamma}^2(x-y)\xi^2_0(y)h^2_0(y)dy\right)^{1/2}.
$$
\end{lemma}

We first discuss the case when $Z^k_t$ are independent one-dimensional Wiener processes.
\begin{thm}
    \label{thm last1}
   Let $Z^k_t$ be independent one-dimensional Wiener processes. Suppose  (\ref{eqn 06.04.1}) and Assumption \ref{ass 06.06.1} hold.
                  Also assume $\gamma\in (-\alpha, \frac{-1-\alpha}{2})$,  $u_0\in U^{\gamma+\alpha-\alpha/p}_p$ and
\begin{equation}
I(p, T):=\left(\E\int^T_0\left(\|f(t,\cdot,0)\|^p_{H^{\gamma}_p}+\|\bar{h}(t,\cdot,0)\|^p_{p}\right)ds \right)^{1/p}<\infty, \label{e:asasqw}
\end{equation}
where
$$
\bar{h}(t,x,0):=\left(\int_{\bR}R_{\gamma}^2(x-y)\xi^2(y)h^2(t,y,0)dy\right)^{1/2}.$$
Then equation (\ref{space-time}) with initial data $u_0$ has a unique solution $u\in \cH^{\gamma+\alpha}_p(T)$ and for this solution,
$$
\|u\|_{\cH^{\gamma+\alpha}_p(T)}\leq c\left(I(p, T)+\|u_0\|_{U^{\gamma+\alpha-\alpha/p}_p}\right).
$$
\end{thm}
\pf
We check whether $f(u)$ and $g(u)$ satisfy condition (\ref{eqn 03.18.91}). Since $\gamma<0$ and $\gamma+\alpha>0$,
$$
\|f(u)-f(v)\|_{H^{\gamma}_p}\leq \|f(u)-f(v)\|_{L_p}\leq K\|u-v\|_{L_p}\leq \varepsilon \|u-v\|_{H^{\gamma+\alpha}_p}+K(\varepsilon)\|u-v\|_{H^{\gamma}_p}.
$$
Also for $g(u)=\{g^k(u)\}$, by Lemma \ref{l last1},
$$
\|g(0)\|_{H^{\gamma+\alpha/2}_{p}(\ell_2)}\leq \|R_{\gamma}\|_{2r}\|\xi\|_{2s}\|h(0)\|_{p}\leq c\|h(0)\|_{p},
$$
$$
\|g(u)-g(v)\|_{H^{\gamma+\alpha/2}_p(\ell_2)}\leq \|R_{\gamma}\|_{2r}\xi\|_{2s}\|h(u)-h(v)\|_{p}\leq c\|u-v\|_{L_p}\leq \varepsilon \|u-v\|_{H^{\gamma+\alpha}_p}+K(\varepsilon)\|u-v\|_{H^{\gamma}_p}.
$$
Therefore condition (\ref{eqn 03.18.91}) is satisfied and the theorem is proved.
\qed

Now we consider space-time white noise with jump L\'evy processes.
Unlike Theorem \ref{thm last1}, in the case space-time white noise with jump L\'evy processes,  $L_p$-theory is not satisfactory due to the condition
$\varepsilon_1>\alpha(1/2-1/p)$ if $p>2$.
 Thus we only give an $L_2$-theory.
\begin{thm}
                  \label{thm last2}
 Suppose $Z^k_t$ are independent one-dimensional jump L\'evy processes with  L\'evy measure $\nu$.  Suppose
 (\ref{main assump}), (\ref{eqn 06.04.1}) and Assumption \ref{ass 06.06.1} hold with $p=2$.
  Also assume $\gamma\in (-\alpha, \frac{-1-\alpha}{2})$,  $u_0\in U^{\gamma+\alpha-\alpha/2}_2$ and $I(2,T)<\infty$, where $I(2,T)$ is taken from (\ref{e:asasqw}).
Then equation (\ref{space-time}) with initial data $u_0$ has a unique solution $u\in \cH^{\gamma+\alpha}_2(T)$ and for this solution,
$$
\|u\|_{\cH^{\gamma+\alpha}_2(T)}\leq c\left(I(2, T)+\|u_0\|_{U^{\gamma+\alpha-\alpha/2}_2}\right).
$$
\end{thm}
\pf
 There is nothing to prove since conditions on $f$ and $g$ were already checked in the proof of Theorem \ref{thm last1}.
\qed

For a stopping time $\tau$  relative to  $\{\cF_t\}$, denote
$$
(\![ 0, \tau]\!]:=\{(\omega,t): 0<t\leq \tau(\omega)\}.
$$ Then obviously the process
 ${\bf 1}_{(\![ 0, \tau]\!]}(\omega,t)$  is left-continuous  and  predictable.  For an $H^{\gamma}_p$-valued
$\cP^{dP\times dt}$-measurable process $u$, write $u\in \bH^{\gamma}_p(\tau)$
if
$$\|u\|^2_{\bH^{\gamma}_p(\tau)}:=\E \left[\int^{\tau}_0 \|u\|^2_{H^{\gamma}_p}ds\right]
<\infty.
$$
We define the Banach spaces $\bL_p(\tau)$, $\bL_p(\tau,\ell_2)$ and $\cH^{\gamma}_p(\tau)$ similarly.
The following theorem plays the key role when we weaken condition
\eqref{main assump}
 later in the next section.

\begin{thm}      \label{remark 12.10}
  Let $\tau\leq T$ be a stopping time.
 Fix a constant $\varepsilon_1$ so that $\varepsilon_1=0$ if $p=2$, and $\varepsilon_1>\alpha(1/2-1/p)$ if $p>2$. Then, under
 Assumption \ref{ass 1} and \eqref{main assump},
 for any $f\in \bH^{\gamma}_p(\tau)$,
 $h\in \bH^{\gamma+\alpha/2}_p(\tau,\ell_2)$,
 $g^{\cdot, j}\in \bH^{\gamma+\alpha/2+\varepsilon_1}_p(\tau,\ell_2), 1 \le j \le m$ and $u_0\in U^{\gamma+\alpha/2-\alpha/p}_p$,  equation (\ref{eqn main levy}) has a unique solution $u$ in $\cH^{\gamma+\alpha}_p(\tau)$, and for this solution
\begin{equation}
                     \label{eqn linear levy2}
\|u\|_{\cH^{\gamma+\alpha}_p(\tau)}\leq c\left(\|f\|_{\bH^{\gamma}_p(\tau)}+\|h\|_{\bH^{\gamma+\alpha/2}_p(\tau,\ell_2)}
+\sum_{j=1}^m\|g^{\cdot, j}\|_{\bH^{\gamma+\alpha/2+\varepsilon_1}_p(\tau,\ell_2)}+\|u_0\|_{U^{\gamma+\alpha-\alpha/p}_p}\right),
\end{equation}
where $c=c(p,T,\delta)$ independent of $\tau$.
 \end{thm}

 \pf
 First we prove the existence and (\ref{eqn linear levy2}). Obviously we have
 $$
 \bar{f}:={\bf 1}_{(\![ 0, \tau]\!]}\,f\in \bH^{\gamma}_p(T),\quad \bar{h}:={\bf 1}_{(\![ 0, \tau]\!]}\,h
 \in \bH^{\gamma+\alpha/2}_p(T,\ell_2), \quad
\bar{g^{\cdot, j}}:={\bf 1}_{(\![ 0, \tau]\!]}\,g^{\cdot, j}    \in \bH^{\gamma+\alpha/2+\varepsilon_1}_p(T,\ell_2).
 $$
Let  $u\in \cH^{\gamma+\alpha}_p(T)$ be the solution of
\eqref{eqn linear}
 with $\bar{f}, \bar{h}$ and $\bar{g}$ instead of $f, h$ and $g$ respectively. Then, since $\tau\leq T$, we have $\|u\|_{\cH^{\gamma+\alpha}_p(\tau)}\leq \|u\|_{\cH^{\gamma+\alpha}_p(T)}$, and by
 Theorem \ref{thm main levy},
\begin{eqnarray*}
\|u\|_{\cH^{\gamma+\alpha}_p(\tau)}&\leq& c \left(\|\bar f\|_{\bH^{\gamma}_p(T)}
  +\|\bar h\|_{\bH^{\gamma+\alpha/2}_p(T,\ell_2)}
  +  \sum_{j=1}^m \|\bar g^{\cdot, j}\|_{\bH^{\gamma+\alpha/2+\varepsilon_1}_p(T,\ell_2)}
  +\|u_0\|_{U^{\gamma+\alpha/2-\alpha/p}_p}\right)\\
 &=& c \left(\| f\|_{\bH^{\gamma}_p(\tau)}
  +\| h\|_{\bH^{\gamma+\alpha/2}_p(\tau,\ell_2)}
  +  \sum_{j=1}^m \| g^{\cdot, j}\|_{\bH^{\gamma+\alpha/2+\varepsilon_1}_p(\tau,\ell_2)}
  +\|u_0\|_{U^{\gamma+\alpha/2-\alpha/p}_p}\right).
 \end{eqnarray*}

Now we prove the uniqueness.  Let $u\in \cH^{\gamma+\alpha}_p(\tau)$ be a solution of equation 
\eqref{eqn linear}. Then obviously,
$$
{\bf 1}_{(\![ 0, \tau]\!]}\cdot(\bD u-a(\omega,t)\Delta^{\alpha/2} u)\in \bH^{\gamma}_p(T), \quad {\bf 1}_{(\![0, \tau ]\!]}\cdot\bS_c u \in \bH^{\gamma+\alpha/2}_p(T,\ell_2),\quad
{\bf 1}_{(\![0, \tau ]\!]}\cdot\bS^{\cdot, j}_d u \in \bH^{\gamma+\alpha/2+\eps_1}_p(T,\ell_2).
$$
  According to
Theorem \ref{thm main levy} we can define $v\in \cH^{\gamma+\alpha}_p(T)$ as the solution of
\begin{eqnarray}
                 \label{e: remark 12.10}
  dv&=&(a(\omega,t)\Delta^{\alpha/2} v+{\bf 1}_{(\![ 0, \tau]\!]}\,(\bD u-a(\omega,t)\Delta^{\alpha/2} u))dt
  +  \sum_{k=1}^\infty 1_{(\![0, \tau ]\!]}\,\bS^k_cu \,dW^k_t\nonumber\\
  &&\quad +\sum_{k=1}^\infty \sum_{j=1}^m 1_{(\![0, \tau ]\!]}\,\bS^{k, j}_d u  \,dZ^k_t, \qquad v(0)=u(0).
\end{eqnarray}
 Then for $t\leq \tau$, $d(u-v)=\Delta^{\alpha/2}(u-v)dt$ and
 $(u-v)(0)=0$.
 Therefore by Theorem \ref{thm deterministic}, we
conclude that $u(t)=v(t)$ for all $t\leq \tau$ a.s..  By replacing $u$ by $v$ for $t\leq \tau$,  from (\ref{e: remark 12.10}) we find that $v$  satisfies
\begin{equation}
                           \label{eqn 5.29.1}
dv=\left(a\Delta^{\alpha/2}v+f{\bf 1}_{(\![ 0, \tau]\!]}\right)dt+ \sum_{k=1}^\infty 1_{(\![0, \tau ]\!]}\, h^k \,dW^k_t
   +\sum_{k=1}^\infty \sum_{j=1}^m 1_{(\![0, \tau ]\!]}\,g^{k,j}  \,dZ^k_t, \quad v(0)=u_0.
\end{equation}
We proved that if  $u\in \cH^{\gamma+\alpha}_p(\tau)$ is a solution of equation
\eqref{eqn linear} then  $u(t)=v(t)$ for  all $t\leq \tau$ a.s.. This proves the uniqueness of solution of equation
\eqref{eqn linear}
in the class $\cH^{\gamma+\alpha}_p(\tau)$ because by Theorem  \ref{thm main levy} $v\in \cH^{\gamma+\alpha}_p(T)$ is the unique solution of equation (\ref{eqn 5.29.1}).
The theorem is proved.
\qed

For a stopping time $\tau \leq T$ and $\gamma \in \bR$, write $u\in \bH^{\gamma}_{p,\loc}(\tau)$ if there exists a
sequence of stopping times $\tau_n \uparrow \infty$ so that $u\in
\bH^{\gamma}_p(\tau \wedge \tau_n)$ for each $n$.

The following is a weakened version of \eqref{main assump}.
 \begin{assumption}
            \label{A3.3}
 There exists an integer $N_0\geq 1$ so that
 $\wh{c}_k<\infty$ for all
 integer $k>N_0$.
\end{assumption}

\begin{defn}\label{D:44}
Let $u_0\in U^{\gamma+\alpha-\alpha/p}_p$,
$f(0)\in \bH^{\gamma}_p(\tau)$,
 $ h (0)\in  \bH^{\gamma+\alpha/2}_p(\tau,\ell_2)$ and
  $ g^{\cdot,j}(0) \in  \bH^{\gamma+\alpha/2+\eps_1}_p(\tau,\ell_2)$, $1 \le j \le m$. We say that  $u\in \cH^{\gamma+\alpha}_{p,\loc}(\tau)$ is  a path-wise solution
to \eqref{eqn main levy}
if the followings hold;

\noindent
{\rm (i)} $u\in \bH^{\gamma+\alpha}_{p,\loc}(\tau)$ and $u(t)$  is  right continuous with left limits in $H^{\gamma}_p$ for $t<\tau$ ($a.s.$),

\noindent
{\rm (ii)} for any $\phi\in C^\infty_0 (\bR^d)$, the equality
\begin{align}
(u(t, \cdot ),\phi)=&(u_0,\phi)+\int^t_0a(\omega,s)(u(s, \cdot), \Delta^{\alpha/2} \phi)ds + \int^t_0(f(s , \cdot),\phi)ds
  \nonumber \\
&+
\sum_{k=1}^{\infty}\int^t_0(h^{k}(s, \cdot),\phi)dW^{k}_s +
\sum_{k=1}^{\infty}\sum_{j=1}^m\int^t_0(g^{k,j}(s, \cdot),\phi)dY^{k,j}_s
             \label{e:3.5}
\end{align}
holds for all $t < \tau$ $a.s.$.
\end{defn}
\begin{thm}
                   \label{thm extension2}
Let $\tau\leq T$. Suppose that  Assumptions \ref{ass 1} and  \ref{A3.3} hold.
 Then for any
  $u_0\in U^{\gamma+\alpha-\alpha/p}_p$,
$f\in \bH^{\gamma}_p(\tau)$,
 $ h \in  \bH^{\gamma+\alpha/2}_p(\tau,\ell_2)$,
  $ g^{\cdot,j} \in  \bH^{\gamma+\alpha/2+\eps_1}_p(\tau,\ell_2), 1 \le j \le m$,  there exists a unique
   path-wise solution $u\in \cH^{\gamma+\alpha}_{p,\loc}(\tau)$
to \eqref{eqn main levy}.
In particular, if $\gamma + \alpha > d/p$, then the unique path-wise solution $u$ is
 $C^{\gamma + \alpha-d/p}$-valued process (for $t\leq \tau$) a.s..
\end{thm}
\pf {\bf{Step 1}}.
 First, additionally assume that  \eqref{main assump} holds.
 Then  the existence of path-wise solution under \eqref{main assump} in $\cH_p^{\gamma+\alpha}(\tau)$ (hence in $\cH^{\gamma+\alpha}_{p,\loc}(\tau))$)  follows from Theorem \ref{remark 12.10}.
 Now we show
that the pathwise solution is unique in $\cH^{\gamma+\alpha}_{p,\loc}(\tau)$.
Let $u\in \cH^{\gamma+\alpha}_{p,\loc}(\tau)$ be a path-wise solution.  Define
$\tau_n=\tau \wedge \inf\{t:\int^t_0\|u\|^2_{H^{\gamma+\alpha}_p}ds>n\}$.  Then $u\in
\bH^{\gamma+\alpha}_p(\tau_n)$ and $\tau_n \uparrow \tau$ since
$\int^t_0\|u\|^2_{H^{\gamma+\alpha}_p}ds<\infty$ for all $t <\tau$, a.s.  By
Theorem \ref{remark 12.10},
$$
\|u\|_{\cH^{\gamma+\alpha}_p(\tau_n)}\leq
c(T,d,\alpha)\left(\|f\|_{\bH^{\gamma}_p(\tau_n)}+\|h\|_{\bH^{\gamma+\alpha/2}_p(\tau_n,\ell_2)}
+\sum_{j=1}^m\|g^{\cdot, j}\|_{\bH^{\gamma+\alpha/2+\varepsilon_1}_p(\tau_n,\ell_2)}+\|u_0\|_{U^{\gamma+\alpha-\alpha/p}_p}\right).
$$
By letting $n \to \infty$ we find that $u\in \cH^{\gamma+\alpha}_p(\tau)$, and the
uniqueness of the pathwise solution under  \eqref{main assump}
follows from the uniqueness result of Theorem \ref{remark 12.10}.

\medskip

{\bf{Step 2}}. For the general case,  note that for each $n>0$ and $ k\leq N_0$,
$$ \wh c_{k,n}:= \left(\int_{\{z\in \bR^m: |z|\leq n\}} |z|^2 \nu_k (dz) \right)^{1/2} \vee \left(\int_{\{z\in \bR^m: |z|\leq n\}} |z|^p \nu_k (dz) \right)^{1/p} <\infty.
$$
Consider   L\'evy
processes $(Z^1_n, \cdots,  Z^{N_0}_n, Z^{N_0+1},\cdots)$ in place of $(Z^1, Z^2\cdots)$, where $Z^k_n (k\leq N_0)$ is  obtained from $Z^k$ by removing all the jumps that has
absolute size strictly large than $n$. Note that condition (\ref{main assump}) is valid with $\wh {c}_k$  replaced by $\wh {c}_{k,n}$.
By {\bf Step 1},
 there is a
unique path-wise solution $v_n\in  \cH^{\gamma+\alpha}_p(\tau)$
 with  $Z^k_n$ in place of $Z^k$ for
$k=1, 2, \cdots, N_0$. Let $T_n$ be the first time that one of the
L\'evy processes $\{Z^k, 1\leq k\leq N_0\}$ has a jump of (absolute)
size in $(n, \infty)$. Define $u(t)=v_n(t)$ for $t <T_n\wedge \tau$. Note that for $n<m$, by {\bf Step 1}, we have $v_n(t)=v_m(t)$  for $t< T_n\wedge \tau$. This is because, for $t<T_n\wedge \tau$,  both $v_n$ and $v_m$
satisfy (\ref{e:3.5}) with each term inside the stochastic integral multiplied by $1_{s<T_n}$ (and with $Z^k_n$, $k\leq N_0$, in place of $Z^k$).
Thus $u$ is well defined.
By letting $n\to \infty$, one constructs a  unique  pathwise solution $u$ in
$\cH^{\gamma+\alpha}_{p,\text{loc}}(\tau)$.
The last claim follows from Sobolev embedding theorem.
The theorem is proved.
\qed

\end{document}